\documentclass[a4paper]{article}
\usepackage[affil-it]{authblk}
\usepackage{epsfig,epic,eepic,psfrag,graphicx,amssymb}
\newcommand{\C}{\ensuremath{\mathbb{C}}}

\newcommand{\HHH}{\mathbb{H}^3}
\newcommand{\HH}{\mathbb{H}^2}
\newcommand{\R}{\mathbb{R}}
\newcommand{\RP}{\mathbb{R}}
\newcommand{\uc}{\ensuremath{\tilde}}

\newcommand{\hs}{\hskip 1 mm}

\newcommand{\chat}{\ensuremath{\mathbb{\widehat{\mathbb{C}}}}}
\newcommand{\qed}{\hfill \ensuremath{\Box}}

\newtheorem{thm}{Theorem}[section]

\newtheorem{lem}[thm]{Lemma}
\newtheorem{prop}[thm]{Proposition}
\newtheorem{cor}[thm]{Corollary}
\newtheorem{definition}[thm]{Definition}
\newtheorem{notation}[thm]{Notation}
\begin{document}

\title{Some Graftings of Complex Projective Structures with Schottky Holonomy \thanks{Partially funded through grant HM1582-08-1-0041 .}
}
\author{Joshua Thompson}

\affil{Department of Mathematics and Computer Science,\\ Northern Michigan University} 


\maketitle

\begin{abstract}

Let $\mathcal{G}^*(S,\rho)$ be the graph whose vertices are marked complex projective structures with holonomy $\rho$ and whose edges are graftings from one vertex to another.  If $\rho$ is quasi-Fuchsian, a theorem of Goldman implies that $\mathcal{G}^*(S,\rho)$ is connected.  If $\rho(\pi_1(S))$ is a Schottky group Baba has shown that $\mathcal{G}(S,\rho)$ (the corresponding graph for unmarked structures) is connected.  For the case that $\rho(\pi_1(S))$ is a Schottky group, this paper provides formulae for the composition of graftings in a basic setting.  Using these formulae, one can construct an infinite number of (standard) projective structures which can be grafted to a common structure.  Furthermore, one can construct pairs of projective structures which can be connected by grafting in an infinite number of ways.


\end{abstract}


\section{Introduction}
\label{intro}

A (complex) projective structure structure on a closed orientable surface $S$ is a generalization of a hyperbolic structure on $S$. A projective structure defines a \textit{holonomy representation} $\rho$ via a \textit{developing map} $D$ and the pair $(D,\rho)$  uniquely determines the structure up to the action of $PSL_2(\mathcal{C})$ by \[(D,\rho) \rightarrow (\phi \circ D, \phi \rho(\gamma) \phi^{-1}) \mbox{ for } \phi \in PSL_2(\C).\]   It is 
natural to ask to what degree the holonomy representation characterizes a projective structure.  Indeed,  Hejhal showed  \cite{Hej} that the holonomy map $\mathcal{H} : P(S) \rightarrow \chi(S)$, where $\chi(S)$ is the $PSL_2(\C)$-character variety of $\pi_1(S)$ and $P(S)$ is the space of all projective structures on $S$, 
 is a local homeomorphism.  Later Earle \cite{Ear} and Hubbard \cite{Hub} independently showed that $\mathcal{H}$ is holomorphic. Also Gallo, Kapovich and Marden \cite{GKM} showed that almost all representations of $\pi_{1}(S)$ 
 into $PSL_{2}(\mathbb{C})$ are holonomy representations of projective structures.


 With almost every representation of the fundamental group of a 
 surface being a holonomy representation of a projective structure we 
 study to what degree the holonomy determines the projective 
 strucuture.  Even when the representation
is discrete and faithful there are many different 
projective structures with the same holonomy representation. 

Examples are constructed using a cut-and-paste  
operation called \textit{grafting}
which glues a projective structure
on an annulus to a projective structure $\Sigma$ on $S$.  This is 
done by first splitting $S$ open along a suitable curve and attaching 
the annulus to $S$ along its boundary.  There are many
suitable curves for grafting, and this gives rise to many different 
projective structures.  A key 
feature of $(2 \pi-)$ grafting, the case that the height of the cylinder is $2\pi$, is that the surgery does not change the holonomy representation.  Therefore the result of $(2 \pi -)$ grafting $\Sigma$ is again a 
projective structure on $S$ with the same holonomy as $\Sigma$.  We only consider grafting in the $2 \pi$ case.

A {\it standard} projective structure is one whose developing map is a covering map of its image.  Baba \cite{Baba} has given a classification analogous to Goldman's theorem (Theorem \ref{Goldman 1} below) for the case when the holonomy is onto a Schottky group whose rank is equal to that of the surface.  Recall that a {\it marking} for $S$ is a homeomorphism $f:S \rightarrow X$ to a fixed surface $X$ of the same genus as $S$.

\begin{thm} [Baba] \label{Baba} Every (un-marked) projective structure on $S$ with Schottky holonomy $\rho$ is obtained by grafting a standard structure once along a multiloop on $S$.
\end{thm}

This theorem implies that $\mathcal{G}(S,\rho)$ is connected if $\rho(\pi_1(S))$  is a Schottky group whose rank is equal to the genus of $S$.   The lack of uniqueness in Theorem \ref{Baba} is due to the fact that, in this case, there are infinitely many standard projective structures with the same holonomy.  

\begin{figure}[ht]
	\begin{center}
  	
	\psfrag{1}{$S$}
  	\psfrag{2}{$\chat$}

 	\includegraphics[scale=.5]{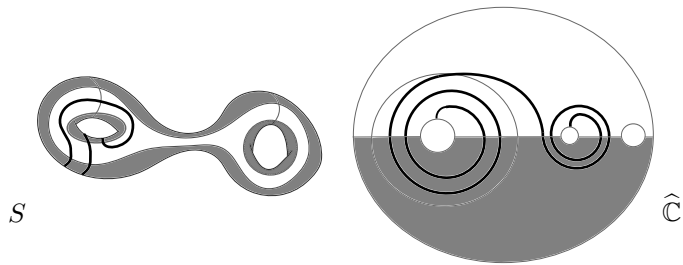}
 
  	\caption{Left-spiraling (in black) curves on a genus two surface and their developed image in $\chat$.  The boundary between grey and white regions are the real curves of the projective structure.}
  
	\label{fig: spiraling_intro}
	\end{center}
  \end{figure}

  Our main result (Theorem \ref{flatsharp}) gives a characterization of grafting along (spherical) spiraling curves (see Figure \ref{fig: spiraling_intro}) in terms of the homotopy type of essential simple closed curves.  Two projective structures differ by an {\it elementary twist} if they differ by a power of a Dehn twist about a curve with trivial holonomy that intersects the {\it real curves} of one of the structures exactly twice. We will show that elementary twists preserve spiraling curves, so the main theorem applies to elementary twists of spiraling curves.  Therefore any two structures differing by an elementary twist may be grafted along (paired) spiraling curves so that the resulting structures are equal.  Furthermore, for such structures there are an infinite number of distinct (paired) spiraling curves suitable for grafting.

\begin{thm} \label{iterated connectivity 1}
 If $\Sigma$ and $\Sigma'$ are standard projective structures with Schottky  holonomy that differ by an 
 elementary twist then there are an infinite number of structures which 
 are grafts of both.
  \end{thm}

The result above implies that $\pi_1(\mathcal{G}^*(S,\rho))$ is not finitely generated.  Since there are an infinite number of standard structures which differ by an elementary twist we also show the following. 

\begin{thm} \label{thm: infinite grafted to one}
There are an infinite number of standard projective structures with Schottky holonomy which can be grafted to a common structure.
\end{thm}

This result together with the theorem of Baba give the following characterization of marked projective structures with Schottky holonomy.

\begin{thm} \label{thm: partial classification}
The graph $\mathcal{G}^*(S,\rho)$ is connected. 

\end{thm}

We remark that the combinatorics of grafting along non-spiraling curves is not explained easily in terms of spiraling curves, therefore a complete classification of marked projective structures with Schottky holonomy in terms of grafting is beyond the scope of this paper.

Since the Schottky holonomy is quasi-conformally conjugate to a real representation from $\pi_1(S)$ onto a Fuchsian group in $PSL_2(\R)$ it suffices to consider the case where the holonomy representation is real.  In this case, $\rho$ preserves $\R$ and thus 
$D^{-1}(\RP)$ is a $\pi_{1}(S)$ invariant set in $\uc{S}$ which
descends to a set of disjoint simple closed curves on 
$S$. So two standard projective structures $\Sigma_{1}$, $\Sigma_{2}$
on $S$ differ by an \textit{elementary 
twist} if the real curves of $\Sigma_{1}$ differ from the real curves 
of $\Sigma_{2}$ by an elementary twist.  These so-called \textit{real curves} characterize a
projective structure with real holonomy in the 
following way.

\begin{thm} \label{real lamination 1} Fix an oriented surface $S$ and a
 real representation $\rho$. 
 Let $\Sigma_{1} = (D_{1}, 
 \rho)$ and $\Sigma_{2} = (D_{2}, \rho)$ be projective structures on 
 $S$.  Then  $\Sigma_{1}$ is isomorphic to $\Sigma_{2}$ if and only if $ D_{1}^{-1}(\R)$ is homotopic to $D_{2}^{-1}(\R)$ on $S$. 
\end{thm}%

 Let $T(S)$ denote the space of all hyperbolic structures on $S$ and $\mathcal{ML}(S)$ the set of measured laminations on $S$.  The representation $\rho$ determines an $X \in T(S)$.  Since $\rho$ is real we assign a measure of $\pi$ to each leaf in $D^{-1}(\R)$.  Theorem \ref{real lamination 1} can be viewed as a special case of the following theorem of Thurston which shows that a projective structure may be viewed as a hyperbolic surface in 
  $\HHH$ which is bent along a geodesic lamination. 

  \begin{thm} [Thurston] There is a canonical homeomorphism $\theta: P(S) \rightarrow T(S) \times \mathcal{ML}(S)$ where $P(S)$ denotes the set of all projective structures on $S$. 
   \end{thm}

One application of Theorem \ref{real lamination 1} is a formulation of
a topological proof of a theorem of Goldman which was mentioned earlier.  This theorem, which motivated much of this work, implies that $\mathcal{G}^*(S,\rho$) is connected if $\rho: \pi_1 \rightarrow PSL_2(\R)$ is injective and discrete. 

\begin{thm}[Goldman]\label{Goldman 1}
Let $\rho$ be a faithful and discrete real representation and
let $\Sigma_{0}$ be the standard projective structure with 
holonomy $\rho$.
For any complex projective structure $\Sigma$ with holonomy $\rho$,
there exists a unique collection of simple closed 
curves $\mathcal{C}$
in $\Sigma_{0}$ such that grafting $\Sigma_{0}$ by $\mathcal{C}$ 
obtains $\Sigma$.
    \end{thm}

Next, we define a special class of admissible curves (see Definition \ref{def: admissible} below) that are the main focus of this paper.

\begin{definition} \label{def: spiraling} 
Choose an admissible curve $\gamma$ and conjugate $\rho$ so that $$Fix(\rho(\gamma) ) = \{0,\infty\}.$$    Let $\uc{\gamma}$ be a component of the preimage of $\gamma$ in the  universal cover of $\gamma$, and set $\hat{\gamma} = D(\uc{\gamma})$.  Then $\gamma$ is spiraling if consecutive points of $\hat{\gamma} \cap \RP$ alternate between $\R^{-}$ and $\R^{+}$.  The action of $\rho(\gamma)$ induces an orientation on the spherical spiral and one of the two fixed points of $\Gamma$ is an attracting fixed point and the other is a repelling fixed point.  If by following the orientation, the spherical spiral winds around the repelling fixed point clock-wise (respectively, counter clock-wise) we say it is \textit{left-spiraling} (respectively, right-spiraling).  The spiraling direction about each fixed point depends only on the homotopy class of $\gamma$. 

\end{definition} 


 An example of a left-spiraling curve is given in Figure \ref{fig: spiraling_intro}.  Theorem \ref{real lamination 1} uses a characterization of grafting along spiraling curves. The real curves obtained by grafting along such curves are 
obtained by a surgery 
procedure, called flat ($\flat$) and sharp ($\sharp$) operations  on the union of the initial real curves and curves isotopic  
to the grafting curve. Such operations are used by Ito in \cite{Ito} to prove the following.

\begin{thm} [Ito] \label{Ito} Let $Q(S)$ be the subset of $P(S)$ of projective structures with quasi-Fuchsian holonomy.  A non-standard component of $Q(S)$ (consisting of quasi-Fuchsian projective structures that are not standard) is not a topological manifold with boundary and any two components of $Q(S)$ have intersecting closures.
\end{thm}

This locally defined operation first removes
intersections and then reconnects nearby curves 
so that the resulting curves are disjoint.  Its relation to 
spiraling curves is expressed below.

  \begin{figure}[ht]
	\begin{center}
	  \psfrag{5}{$Gr(\gamma_l$)}
	  \psfrag{4}{$Gr(\gamma_r$)}
	  \psfrag{1}{$\gamma_l$}
	  \psfrag{2}{$\gamma_r$}
	  \psfrag{3}{$[\lambda,2 \gamma]_\sharp$}
	  \psfrag{6}{$[\lambda,2 \gamma]_\flat$}

 \includegraphics[height=200px,width=330px]{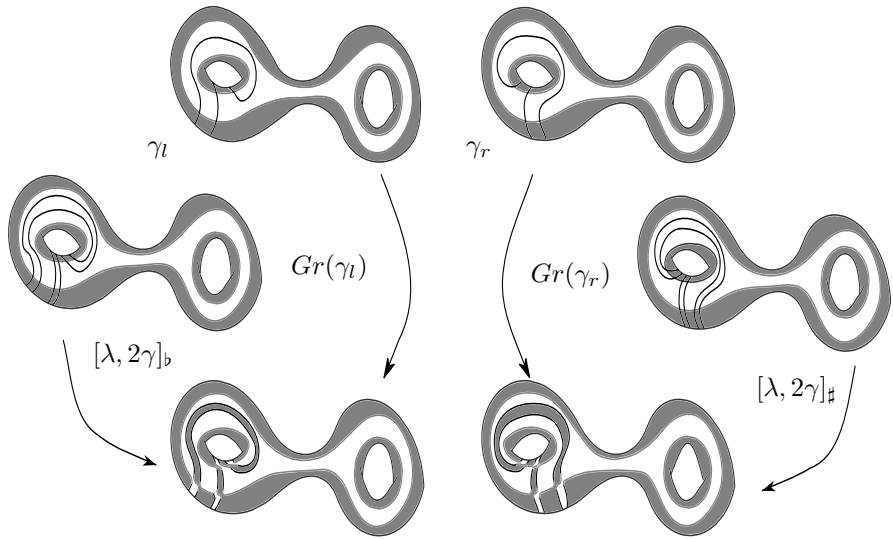}
 
  \caption{The effect of grafting along left and right-alternating curves expressed as  $\flat$ and $\sharp$, respectively, operations on the union of $2 \gamma$ with the real curves of the original structure.}
  \label{fig: sharp_flat_theorem}
	\end{center}
	\end{figure}

\begin{thm} \label{flatsharp 1}
If $\gamma$ is a spiraling admissible curve in $\Sigma$,
then the real curves of the graft of $\Sigma$ along $\gamma$ are 
obtained by either a $\sharp$ or $\flat$ operation on the union of the
real curves of $\Sigma$ 
with two curves isotopic to $\gamma$.  
\end{thm}

\subsection{Idea of Proof of Theorems \ref{iterated connectivity 1} and \ref{thm: infinite grafted to one}}

By Theorem \ref{real lamination 1} complex projective structures with real holonomy are classified by the homotopy classes of their real curves.  Theorem \ref{flatsharp 1} shows that real curves of certain (spiraling) grafted structures are obtained from the original real curves using surgery operations.  We will show that the $\sharp$ and $\flat$ operations on spiraling curves in $S$ reduce to $\sharp$ and $\flat$ operations on curves in a torus.  In the torus these operations are closely related to Dehn twists.  This relationship is used to show the real curves produced by grafting can be constructed by Dehn twists of the inital real curves about a meridian in $S$ that intersects the real curves exaclty twice.  This technical condition on $\beta$ makes it easy to see the relationship to spiraling curves.  Powers of Dehn twists correspond to higher frequencies of spiraling.  The meridian is used to identify other projective structures, via an elementary twist, and we will show that certain grafts of these other structures are also grafts of the initial structures.

For any fixed real lamination $\lambda$ of a standard structure $\Sigma$ and meridian $\beta$ which intersects $\lambda$ exactly twice, we will show that the each power of a Dehn twist of $\lambda$ about $\beta$ are the real curves of some standard structure $\Sigma'$.  We will then use the Dehn twist to produce admissible curves for $\Sigma$ and $\Sigma'$ so that the graftings of both coincide.


%

\subsection*{Acknowledgements}
I thank Ken Bromberg for introducing me to this beautiful area of mathematics, Shinpei Baba for his careful observations, Renzo Cavalieri for help with the introduction, and the reviewer for their very thoughtful comments.

\section{Preliminaries} \label{Prelim}



A \textit{complex projective structure} $\Sigma$ on an
orientable surface $S$ is a maximal atlas  
$(\phi_{i} : U_{i} \rightarrow \chat )_{i\in I}$ where $\phi_{i}$ is 
a homemorphism onto its image and the transition functions
are restrictions of elements in $PSL_{2}(\chat)$, to 
connected components. Since a maximal atlas is unique the set of charts giving a projective 
structure is unique.

A global definition of a projective structure  is an equivalence class of
developing pairs $[(D, \rho)]$ where the \textit{developing map} 
$D : \uc{S} \rightarrow \chat$ 
globalizes the charts in the atlas and the homomorphism
$\rho: \pi_{1}(S) \rightarrow 
PSL_{2}(\chat)$ globalizes the coordinate changes.  
Then $\rho$ is a \textit{holonomy} 
representation and commutes with $D$ \[ D(\gamma x)=\rho(\gamma) D(x) \] for all $\gamma \in \pi_{1}(S)$ and $x \in \uc{S}.$ The developing map $D$ is unique up to pre-composition with a map that is a lift of map of $S$ that is isotopic to the identity and post-composition with PSL$_{2}(\chat)$ and $\rho$ is unique up to
conjugation.  To ease notation we let $(D,\rho)$ denote a projective structure where $D$ is a representative developing map and $\rho$ is a representative holonomy representation of the equivalence class $[(D,\rho)]$.   


A basic example of a complex projetive structure on $S$ is the 
quotient of the upper half plane in $\mathbb{C}$ by a Fuchsian group, 
i.e., a discrete subgroup of $PSL_{2}(\R)$ isomorphic to $\pi_{1}(S)$.
The quotient is both a hyperbolic structure and a 
complex projective structure.  Indeed, since Isom($\HH$) $\cong 
PSL_{2}(\R) \subset PSL_2(\mathbb{C})$ and $\HH \subset
\mathbb{C} $ every hyperbolic structure determines a complex 
projective structure.  More generally, if $\Gamma$ is a Kleinian group,
(a finitely generated discrete subgroup of 
$PSL_{2}(\chat)$) isomorphic to $\pi_{1}(S)$ and $\Omega$ is a 
connected component of the domain of discontinuity of 
$\Gamma$ which is fixed by $\Gamma$ then 
the quotient $\Omega \slash \Gamma$
is a complex projective structure on $S$.  However, as well shall see 
developing maps are not necessarily injective and holonomy 
representations are not necessarily discrete or faithful.  

In the next section we first formulate a precise global definition of a 
projective structure, expressed a developing map and holonomy 
representation.  
After mentioning relevant earlier work in this area,
we describe the grafting operation in detail, a 
construction which is independent of the underlying
holonomy representation.  Restricting our discussion to 
projective structures whose holonomy representations are real, we 
then discuss the effect of grafting such structures.

\subsection{The Developing Map}
\label{The Developing Map}

We first give a precise definition of the developing map.
Let $\mathcal{A} = (\phi_{i}, U_{i})$ be a maximal atlas 
defining a projective structure $\Sigma$ on an oriented surface $S$. 
Choose a point $x$ in $S$ and assume $x \in U_{0}$.  A lift $\uc{x}$ 
of $x$ is contained in a lift of $\uc{U}_{0}$ of $U_{0}$ and a chart 
$(\phi_{0},U_{0})$ on $S$ lifts to a chart 
$(\uc{\phi}_{0},\uc{U}_{0})$ on $\uc{S}$. 
Since $\uc{S}$ is simply 
connected $(\uc{\phi}_{0},\uc{U}_{0})$ extends to a chart 
$D_{0}$ on all of $\uc{S}$.  This chart is maximal since its domain 
is all of $\uc{S}$.

This global chart is unique in the sense 
that if two global charts agree on an open set of $\uc{S}$ then there 
is an an isotopy of $\uc{S}$, after which they agree on all of 
$\uc{S}$.  Therefore, this global chart is unique up to an isotopy of 
$\uc{S}$.

For another lift $\uc{x}_{g}$ of $x$, there is some $g$ in 
$\pi_{1}(S)$ such that $g \uc{x} = \uc{x}_{g}$ and $g \uc{U}_{0} = 
\uc{U}_{g}$ is 
a lift of $U_{0}$.  There is a lift $(\uc{\phi}_{g},\uc{U}_{g})$ 
of $(\phi_{0},U_{g})$ which extends to a maximal global chart $D_{g}$ on 
$\uc{S}$.  Since $D_{0}$ and $D_{g}$ arise from lifts of $U_{0}$ they 
differ by pre-composing with an element of the fundamental group and 
so there is some $g$ such that \[ D_{0} = D_{g} \circ g .\] Both $D_{g}$ and $D_{0}$ are global charts on $\uc{S}$ and so there is a 
transition map $\phi_{g}$ which we post-compose with $D_{0}$ so that 
the composition agrees with $D_{g}$ on all of $\uc{S}$, i.e., \[  D_{0} = \phi_{g} \circ D_{g} .\]  We call each $D_{i}$ a \textit{developing map}.  
This determines a map $\rho : \pi_{1}(S) \rightarrow PSL_{2}(\C)$ 
given by 
\[ \rho (g) = \phi_{g} .\]
We show that $\rho$ is a homomorphism.  If $g_{1}$ and $g_{2}$ 
are elements of $\pi_{1}(S)$, then $D_{0} \circ g_{1} \circ g_{2}$ is 
a global chart on $\uc{S}$.  By the definition of an atlas there is a
Mobius transformation $\phi$
such that 
\[ \phi \circ D_{0} = D_{0} \circ g_{1} \circ g_{2}
.\]
Similarly, there are Mobius transformations $\phi_{1}$ and 
$\phi_{2}$ such that 
\[  \phi \circ D_{0} = \phi_{2} \circ 
D_{0} \circ g_{1} = \phi_{1} \circ \phi_{2} \circ D_{0} .\]
Therefore $\phi = \phi_{1} \circ \phi_{2}$ and thus $\rho$ is a 
homomorphism we call the \textit{holonomy representation}.

If another point $y$ in $S$ is chosen, the  resulting developing map
differs from $D_{0}$ by a Mobius transformation since both are global 
charts on $\uc{S}$.  Likewise the resulting holonomy representation 
differs from the initial one by conjugation in $PSL_{2}(\C)$. 

 Fix a projective structure $\Sigma(D, \rho)$ on $S$.
 If $X$ is a subset of $S$, and $\uc{X}$ is a lift of $X$ to the 
 universal cover, then we say $\hat{X} = D(\uc{X})$ is 
 the developed image 
 of $X$.


\subsection{Grafting}
   
The process of grafting is defined for an admissable curve
$\gamma$ in $S$, and a projective structure $\Sigma = (D, \rho)$.  

\begin{definition}
\label{def: admissible}
 A simple closed curve $\gamma$ in $S$ is admissible if :
 \begin{itemize}
  \item $\rho(\gamma)$ is hyperbolic and 
  \item $D \arrowvert_{\uc{\gamma}}$ is an embedding of each lift 
  $\uc{\gamma}$ of $\gamma$
  \end{itemize}
  where $\rho(\gamma)$ denotes the image under $\rho$ 
  of an element of $\pi_{1}(S)$ in the free homotopy class of $\gamma$.
  We say the homotopy class $[\gamma]$ of $\gamma$ has an admissible 
  representative if there is an admissible curve in $[\gamma]$.  
 \end{definition}

\subsection{Hopf Torus}

To every grafting there is an associated torus, which we now briefly discuss, and discuss in more detail later.  Since $\gamma$ is admissible  $\rho(\gamma)$ is a hyperbolic element 
in $PSL_{2}(\mathbb{C})$.  There are two fixed points 
$Fix(\rho(\gamma))$ of $<\rho(\gamma)>$.
The quotient $T = \chat \backslash Fix(\rho(\gamma))$
under the action of $<\rho(\gamma)>$ is a projective structure on 
a torus we call the Hopf torus.  
Each lift of $\gamma$ has a corresponding Hopf torus. The developing 
maps of these Hopf tori differ by post-composition with a Mobius transformation and their 
holonomy representations differ by conjugation in $PSL_{2}(\C)$.

\subsection{The Collapsing Map} \label{collapse}

\begin{figure}[ht]
\centerline{
  \epsfig{file=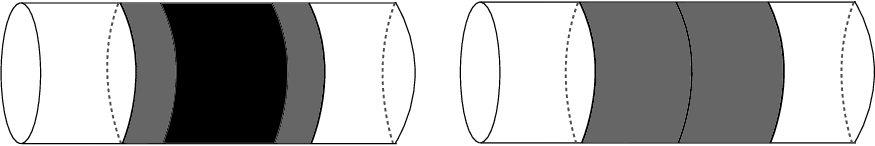,scale = .5}
}
\caption{The collapsing map}
\label{figure: collapsing map}
\end{figure}

Let $N(\gamma)$ denote an annular
neighborhood of an admissible curve $\gamma
\in S$ and let $\eta : S^{1} \times (0,1) \rightarrow N(\gamma)$ be 
a homeomorphism such that $\eta(S^{1}, \frac{1}{2}) = \gamma$.  We define a
\textit{collapsing} map $\nu : S \rightarrow S$ which is the identity
on $S \backslash N(\gamma)$, collapses $\eta (\gamma \times (\frac{1}{4},\frac{3}{4}))$ to 
$\gamma$ and expands 
$\eta(\gamma \times (0,\frac{1}{4}))$ and $\eta(\gamma \times (\frac{3}{4},1))$ as follows and as shown in Figure \ref{figure: collapsing map}.

 \[ \nu (y) =
\left \{ \begin{array}{ll}
	y & \mbox{if $y \in S \backslash N(\gamma)$} \\
	\eta(s, \frac{1}{2}) & \mbox{if $y = \eta(s , t) $ for 
	$t \in (\frac{1}{4},\frac{3}{4})$ } \\
	\eta(s,2t) & %
	  \mbox{if $ y = \eta(s , t)$ for %
	    $t \in (0,\frac{1}{4})$ } \\
	\eta(s,2t-1) & %
	   \mbox{if $ y = \eta(s , t)$ for %
	$t \in (\frac{3}{4},1) $ } 
    \end{array}
	\right . \]

 Choose a lift to the universal cover $\uc{\nu} : \uc{S} \rightarrow
\uc{S}$ which is homotopic to the identity outside of the preimage 
of $N(\gamma)$.

\subsection{Projective Structure on Annulus} \label{annulus 
structure}

Next, we define a projective structure on an annulus which is used in grafting to extend the projective structure on the base surface to the annulus.  We define the developing map on an annulus $A$ by choosing a lifting a particular map $h: A \rightarrow T$ from $A$ to a torus $T$ which identifies the boundary of $A$.  We extend the lift equivariantly and show that the resulting developing map does not depend on the choice of lift. 

Choose a lift $\uc{\gamma}$ of $\gamma$ in $\uc{S}$.  
Let $A = \eta(\gamma \times (\frac{1}{4},\frac{3}{4}))$ and $\uc{A}$
be the lift to $\uc{S}$ containing $\uc{\gamma}$.
Since
$\gamma$ is admissible $D(\uc{\gamma})$ is a
simple arc in $\chat$ with distinct endpoints $p_1$ and $p_2$ .  For $x \in S^1$ let $a_x = \eta (x \times (\frac{1}{4},\frac{3}{4}))$.  
There is a unique non-trivial homotopy class of simple closed curves in $\chat \slash Fix(\rho(\gamma))$, and we call this class $[c]$. 
Let $h$ denote an orientation 
preserving map to the torus $T$ and let
\[h: A \rightarrow T\] which has a lift $\uc{h} : \uc{A} \rightarrow 
\chat \slash Fix(\rho(\gamma))$ which is surjective, 
injective away from $\uc{\gamma}$, on each 
component of $\partial(\uc{A})$ restricts to
\[ \uc{h} = D \circ \uc{\nu} \] 
and for each $x \in S^1$ has the property that
\[ \uc{h} ( a_x ) = \delta_x  \] where $\delta_x \in [c]$.  These properties ensure that arcs of $A$ that are collapsed to single points of $\gamma$ develop to simple closed curves in $\chat \slash Fix(\rho(\gamma))$.  Since $\uc{h}$ is a lift of $h$, $\uc{h}$ is equivariant 
with respect to the action of $<a>$ the $\mathbb{Z}$-subgroup
that fixes $\uc{\gamma}$.  That is, for every $x \in \uc{\gamma}$ 
\[ \uc{h}(ax) = \rho(a)\uc{h}(x)
.\]
%


Assume $x$ lies in $\pi^{-1}(A)$ where $\pi: \uc{S} \rightarrow S$ is the 
universal covering. Then there is some non-trivial
$g_{1} \in \pi_{1}(S)$ for which $g_{1}x$ lies in $\uc{A}$.
We extend 
$\uc{h}$ to $\pi^{-1}(A)$ by 
\[ \uc{h}(x) = \rho(g^{-1}_{1}) \uc{h} (g_{1}x)
 .\]

Note that there are many $g$ such that $gx \in \uc{A}$.  If $g_{2}x 
\in \uc{A}$ note that $g_{1} g^{-1}_{2}$ fixes $\uc{A}$.  Then 
\[ \rho(g^{-1}_{2}) \uc{h}(g_{2} x) = \rho(g^{-1}_{2}) \rho(g_{2} 
g^{-1}_{1}) \uc{h}(g_{1} g^{-1}_{2} g_{2} x) = \rho(g^{-1}_{1} )
\uc{h}(g_{1}x) .
\]
Therefore the definition of $\uc{h}$ is independent of the choice of 
$g$.  For any $x \in \pi^{-1}(A)$, and $g \in \pi_{1}(S)$ 
we have 
\[ \uc{h}(gx) = \rho(g^{-1}_{1})\uc{h}(g_{1}gx) =
   \rho(g^{-1}_{1}g_{1}) \uc{h}(gx) =
     \rho(g) \uc{h}(x). \]
Then $\uc{h}$ is equivariant with respect to the action of the 
fundamental group of $S$.


\subsection{Grafted Developing Map}

Given a projective structure on $S$ and
the collapsing map $\nu$ defined above we define a new developing map $D'
: \uc{S} \rightarrow \chat$  in the following way.  For $\uc{A}$ 
a lift of $A$ to $\uc{S}$ we set

\[ D' (x) = \left
\{ \begin{array}{ll}
		D \circ \uc{\nu} (x) & \mbox{ if $x \in \uc{S}
		\hs \backslash \hs \pi^{-1}(A)$} \\
		\uc{h}(x) & \mbox{ if  $x \in \pi^{-1}(A)$}
		\end{array}
	\right . \]
	
Since $\nu$ is a homotopy equivalence we have the commutativity
relation $\uc{\nu}(gx) = g\uc{\nu}(x)$ for all $g \in \pi_1(S)$.  
For $g \in \pi_1(S)$, $gx$ lies in $\pi^{-1}(A)$ if and only if $x$ 
lies in $\pi^{-1}(A)$ so we have 
\[ D' (gx) = \left \{ \begin{array}{ll}
	D \circ \uc{\nu} (gx) = D (g \uc{\nu} (x)) = \rho(g)
	D'(x) & \mbox{ if $gx \in \uc{S} \hs \backslash \hs
	\pi^{-1}(A) $}\\
	\uc{h} (gx) = \rho(g) \uc{h}(x) = \rho(g) D'(x) & \mbox{ if  $gx \in \pi^{-1}(A) $} 
		\end{array}
	\right . \] which shows that $D'$ is $\pi_{1}$-
	equivariant with holonomy $\rho$.  
 This defines a new projective
structure $(D', \rho) = Gr_\gamma(\Sigma)$ on $S$.


%
%
%

%


\section{Real Holonomy}  \label{real holonomy 2}

This section consists of facts about complex projective 
structures with real holonomy.  In the first section no injectivity 
assumption is made on the holonomy representation.  Here we reduce 
the classification of projective structures with real holonomy to the 
classification of the homotopy class of their 
{\it real curves}.  
By the {\it real curves} of $(D, \rho)$ we refer to all 
curves in the surface whose lift to the universal cover lies in
the preimage of
$\RP$ under $D$.  This reduction is used to construct an alternate proof of a 
Theorem of Goldman (see \cite{Gol}).


 \subsection{Real Curves of a Projective Structure} \label{real projective}

 In the remainder of this section, we will assume the representation $\rho$
is real.  To ease notation we make the following convention.
 
 \begin{notation}
    Let $\Sigma = (D,\rho)$ be a projective structure on $S$ and 
    let $\pi: \uc{S} 
    \rightarrow S$ be the universal covering projection.  Set
    \[ \lambda = \pi( D^{-1}(\R) ) \] and let $\uc{\Sigma}$ denote 
    the universal cover of the projective structure $\Sigma$.
    \end{notation}

    Next, we note that the real curves of a projective structure are closed curves 
 on the surface.  This was used implicitly in \cite{Tan88}.

 \begin{lem} \label{closed curves}
  Let $(D,\rho)$ denote a projective structure on $S$, and $\pi: 
  \uc{S} \rightarrow S$ the universal covering.  Then 
  $\lambda$ is a finite collection of disjoint closed curves 
  on $S$.
  \end{lem}
 
  \noindent \textbf{Proof : }  
  
  Since $D$ is continuous and $ \R \cup \{\infty\}$ is closed in \chat, 
  it follows that $\uc{\lambda} = D^{-1}(\R \cup \{\infty\})$ is closed in $\uc{S}$.  
  The real representation $\rho$ preserves $\R \cup \{\infty\}$ and so 
  $\uc{\lambda}$ is a closed $\pi_{1}$-invariant set in $\uc{S}$.  
  Then $\lambda$ is a closed set in the compact surface $S$, and is 
  therefore compact.  Each point in $\lambda$ has a 
  neighborhood in $S$ homeomorphic to the unit disk in $\C$ with the 
  real line corresponding to the real curve since $D$ is a local 
  homeomorphism.  Therefore each 
  connected component of the real curves is a 1-manifold.  A compact 
  1-manifold is homeomorphic to a circle so each connected component 
  of $\lambda$ is a closed curve.  Since $\lambda$ is compact the number of curves in $\lambda$ is finite.  \qed
  
	A {\it geometric disk} of $\chat$ is an open subset of $\chat$ 
   whose boundary is a geometric circle.  A 
   geometric disk of $\uc{S}$ is an open subset $U$ such that $D: U 
   \rightarrow D(U)$ is a homeomorphism onto a geometric disk of 
   $\chat$.  The set of geometric disks of $\uc{S}$ is partially 
   ordered by inclusions.  We call a maximal geometric disk a maximal 
   disk.

  \begin{lem} \label{unique disk}
   Every 
   point of $\uc{S}$ lies in a proper maximal disk.
   \end{lem}

  Let $D^{2} = \HH \cup S^{1}$ be the compactification of 
  $\uc{\mathbb{H}}^{2}$.  If $U$ is a maximal disk of $\uc{S}$ we call
  \[U_{\infty} = \overline{U} - \uc{S} \]
  the ideal set of $U$ where $\overline{U}$ is 
  the closure of the 
  compactification of $U$ in $\uc{S}$.  Since $\overline{U}$
   is conformally equivalent to a closed ball we may form 
  $C(U_{\infty})$ the convex hull of $U_{\infty}$ in $U$.  Then 
  $C(U_{\infty})$ is the smallest convex set in $\overline{U}$ 
  containing $U_{\infty}$.    
  
	\qed   

  \begin{lem} \label{pleated}
   Every point of $\uc{S}$ lies the convex hull $C(U_{\infty})$ of a 
   unique maximal disk and $\partial(U_{\infty})$ is closed in $\uc{S}$.
   \end{lem}
   
   Let $D^{3} = \HHH \cup S^{2}$ be the compactification of $\HHH$, and let 
   $\mathcal{H}(D(U))$ denote the hyperbolic plane in $\HHH$ that
   contains the  boundary of $D(U)$ for a maximal disk $U$.  
   For each $U$ there is a nearest point retraction $\Phi_{U}: D(U) 
   \rightarrow \HH \subset \HHH$ given by 
   \[ \Phi_{U} : D(U) \rightarrow \mathcal{H}(D(U)) 
   .\]

   By  Lemma \ref{unique disk}
   for  every $x \in C(U_{\infty})$ there is a map
   $\Psi: \uc{S} \rightarrow  \HHH$ given by 
   \[ \Psi(x) = \Phi_{U}(D(x)) .\]
   In general the map $\Psi$ is not
   locally injective.  Denote by $\mathcal{B}$ the set of those 
   $C(U^{i}_{\infty})$ 
   on which $\Psi$ is not injective. 
   There is a lift of a homotopy equivalence $\uc{\nu}$ 
   of $\uc{S}$ which collapses each 
   component of $\mathcal{B}$ to a common arc.  Then 
   \[ \Psi : \uc{\nu}(S) \rightarrow \HH \subset \HHH \]
   is injective and pulls back a hyperbolic metric on $\uc{\nu}(S)$.  Since 
   $\uc{\nu}$ is a homotopy equivalence this induces a hyperbolic metric 
   $\sigma$ on 
   $S$.
    
     It can be seen that $\Psi$ is a {\it pleated} map, that is a continous 
     map from a hyperbolic surface $S$ into a hyperbolic 3-manifold 
     ($\HHH$ in this case) such that for any point $x \in S$, there 
     is a geodesic in $S$ containing $x$ which is mapped to a 
     geodesic in the 3-manifold, and the path metric induced 
     from $\HHH$ by 
     $\Psi$ coincides with the hyperbolic metric on $S$.  The image 
     of $\Psi$ is a pleated surface in $\HHH$ 
     which is bent along a lamination $\lambda$ (the pleating
     locus) in the sense that 
     $\Psi(S) \arrowvert \lambda$ is totally geodesic.  
     The geodesic laminaiton $\lambda$ supports a measure, 
     if $\Psi$ is locally convex and this measure corresponds
     to the bending angle. 
    
\qed

\subsection{The Measured Lamination} \label{the measured lam}

Here we realize	 the real curves as a measured lamination with isolated leaves, each of weight $m \pi$ for $m \in \mathbb{Z}$. Since we assume $\rho$ is real, the limit set of $\Gamma$ is contained in
     $\R \cup \{\infty\}$, and the convex hull $C$ of the limit set 
     $\Lambda$ in $\HHH$ is contained in the equitorial half plane
     $\textbf{H}$.  The \textit{pleating locus} of $\uc{S}$ is the
     set of points $L$ in $\uc{S}$ through which there passes one and
     only one open geodesic arc $a$ for which $\Psi(a)$ is a geodesic
     in $\Psi(\uc{S})$.  It is shown in \cite{MT} p. 177,  that the image
     of the pleated map, a pleated surface $S'$, is contained in the
     convex hull $\mathcal{CH}(\Gamma)$ of $\Gamma$  and that $L$ is a
measured geodesic lamination
on $\uc{S}$.  


Then, following \cite{Pap}, $L$ assigns to any smooth compact arc $c$ on $\uc{S}$ intersecting $L$ transversely, and whose endpoints lie in the complement of $L$ a finite Borel measure $\mu$ on $c$ with support in $c \cap L$.  By Lemma \ref{unique disk} each endpoint $e_i$ of $c$ is containted in the convex hull $C(U_{e_i})$ of unique maximal disk $U_{e_i}$.  Since the $e_i$ do not lie in the pleating locus there are at least two geodesics of $U_{e_i}$ through $e_i$ which are mapped to geodesics by $\Psi$.  Now each $U_{e_i}$ has more than two ideal points and the $C(U_{e_i})$ are two-dimensional.  It follows that $\Psi(C(U_{e_i}) \subset \textbf{H}$ where $\textbf{H}$ denotes the equitorial half plane in $\HHH$.

Let $\Theta(s,t)$ denote the dihedral angle of intersection of the circles $U_s$ and $U_t$.   In \cite{KT} it is shown that $\mu(c) = \varphi'$ with \[ \varphi = \lim \Theta(0,t_1) + \dots + \Theta(t_n,t) \]  where the sum runs over subdivisions $t_1 < t_2 < \cdots < t_n$ of $[0,t] = c$ and the limit is taken as the width of the widest interval goes to zero.

If $\mu(c) < \pi$ then $\Theta(D(U_{e_1}),D(U_{e_2})) \leq \lim (\Theta(0,t_1) + \dots + \Theta(t_n,t)) < \pi$.  Then $\mu(c) = 0$, since otherwise $\Psi(C(U_{e_i}) \subset \textbf{H}$ implies a contradiction.  It follows easily that if $\mu(c) = m \pi$ for $m \in \mathbb{Z}^+$ then $c$ is an atom of $\mu$ and all atoms have measure $m \pi$ and correspond to seperated leaves of the lamination. If $\mu(c) > \pi$ and has no atoms then there is a $c' \subset c$ such that $\mu(c') < \pi$ and we conclude that the support for $\mu$ consists of isolated leaves each with weight $m \pi$.  
   
The proof of the following proposition is technical and the basic idea is summarized as follows.  Fix a connected component  $U \in \uc{S} \backslash D^{-1}(\R)$ and let $x \in U$ be arbitrary. Either $U$ is the unique maximal disk whose convex hull contains $x$, or there is a point $p \in U$ near to $x$ where $U$ is the unique maximal disk whose convex hull contains $p$.  Since the representation is real, maximal disks are mapped onto either the upper or lower half plane.

  \begin{prop} \label{max disks}
	  
   Fix a projective structure $(D,\rho)$ on $S$ where $\rho$ is a 
   discrete real representation onto $\Gamma$.  
   Let  $U$ denote a connected component of 
   $U \subset \uc{S} \backslash D^{-1}(\R)$.  Then $U$ is a 
   maximal disk of 
   $\uc{S}$ which is mapped either onto the upper or lower half plane by $D$.
   \end{prop}

   \noindent \textbf{Proof : }   

   Fix $U$ a connected component of $\uc{S} \backslash D^{-1}(\R)$.  
 By Lemma \ref{unique disk}, for $x \in \uc{S}$ there is a unique maximal
disk $V_x$ such that $x$ lies in $C(V_x)$.  By construction, $D(V_x)$ is a geometric disk in $\chat$.  The complete hyperbolic metric on $D(V_x)$ pulls back to a complete hyperbolic metric on $V_x$.  Let $H_x$ denote the plane in
$\HHH$ whose boundary is the same as the
boundary of $D(V_x)$.  

Since $V_x$ is maximal, $C(V_x)$ has at least two ideal points which bound a
geodesic $g_x$ in the hyperbolic metric on $V_x$ containing $x$.  As argued above, $\rho$ being real implies that  $\Psi(C(V_x)) \subset \textbf{H}$.   The limit points of $V_x$ are mapped into $\R$ by $D$, and thus $H_x$ intersects $\textbf{H}$.  In fact, since $\Psi(\uc{S}) \subset \textbf{H}$ we have \[\Psi(x) \in H_x \cap \textbf{H}.\]  If $H_x = \textbf{H}$ then $V_x = U$ and the proof is complete, otherwise $H_x \cap \textbf{H}$ is a geodesic.    

Since $x \in U$  we can assume $D(x) \in \Omega_+$. Let $g^\perp$ be the ray based at $\Psi(x)$ that is perpendicular to $\textbf{H}$ and has ideal endpoint $p$ in $\Omega_+$. The ray from $\Psi(x)$ to $D(x)$ will make an acute angle with $g^\perp$ so $p$ will be in $D(V_x)$.


Let $y = D^{-1}(p) \cap V_x$, then $y \in U$ and to complete the proof we show that $V_y = U$. As before let $H_y$ be the plane in $\HHH$ whose boundary is $\partial(D(V_y))$.   Since the plane $H_y$ divides $\HHH$ into two half spaces, one of these half spaces $C$ has $D(V_y)$ as its ideal boundary and contains $p$.  

Since $\Psi(\uc{S}) \subset \textbf{H}$, the geodesic in $\HHH$ with endpoint $p$ that is orthogonal to $H_y$ intersects $H_y$ in $\textbf{H}$.  
Note that the horosphere based at $p$ that is tangent to $\mathbb{H}$ intersects $\mathbb{H}$ at $\Psi(x)$. Therefore any horoball based at $p$ that intersects $\mathbb{H}$ will contain $\Psi(x)$. Since the horoball based at $p$ that is tangent $H_y$ intersects $H_y$ at a point in $\mathbb{H}$ this horoball contains $\Psi(x)$. As this horoball will be contained in $C$ we have that $\Psi(x) \in C$ and the geodesic $\Psi(g_x)$ which contains $\Psi(x)$ will also intersect $C$.  Therefore one of these endpoints must be contained in $D(V_y)$ the ideal boundary of $C$.  

We will show that $D$ is injective on $V_x \cup V_y$. First note that since $D$ is a local homeomorphism it is an open map and $D(V_x \cap V_y)$ is an open subset of $D(V_x) \cap D(V_y)$. The continuity of $D$ implies that it is also a closed subset of $D(V_x) \cap D(V_y)$ and since $V_x \cap V_y \neq \emptyset$ it is also not empty. This implies that $D(V_x \cap V_y) =D(V_x) \cap D(V_y)$.

Since $D$ is injective on both $V_x$ and $V_y$ to show that it is locally injective on their union we need to show that if $q_0 \in V_x \backslash V_y$ and $q_1 \in V_y \backslash V_x$ then $D(q_0) \neq D(q_1)$. We proceed by contradiction. If $D(q_0) = D(q_1)$ then the image is in $D(V_x) \cap D(V_y)$. By the above paragraph there is then a $q_2 \in V_x \cap V_y$ such that $D(q_0) = D(q_2)$, contradicting the local injectivity of $D$ on $V_x$.

Now let $z_\infty$ be a point $D(V_x) \cap D(V_y)$ and let $p_\infty$ be the unique point in $V_y$  such that $D(p_\infty) = z_\infty$. There will then be a sequence of points $z_i$ in $D(V_x) \cap D(V_y)$ such that $z_i \rightarrow z_\infty$. Let $p_i$ be the unique point in $V_x \cup V_y$ such that $D(p_i) = z_i$. Since $D(V_x \cap V_y) = D(V_x) \cap D(V_y)$ we have that $p_i \in V_x \cap V_y$. Since $D$ is a homeomorphism on $V_y$ we also have $p_i \rightarrow p_\infty$.  Therefore $z_\infty$ is not an ideal point of $V_x$ and $V_y$ is a maximal disk contained in $U$ where maximality then implies $V_y = U$.

   \qed

  In the next lemma we show that every essential simple closed curve disjoint from the real curves of a projective structure is admissible.

   \begin{lem} \label{disjoint admissibility}
    Let $\rho$ be a discrete real representation. If $\gamma$ is an essential
    simple closed curve disjoint from 
    $\lambda$ in $\Sigma(\rho, \lambda)$, then 
    $\gamma$ is admissible.
    \end{lem}
    
    \noindent \textbf{Proof : }  
     Since $\gamma$ is essential, any lift of $\gamma$ to $\uc{S}$ has 
     two distinct endpoints on the boundary. By Proposition \ref{max 
     disks} any such lift is contained in a maximal 
     disk since $\gamma$ is disjoint from $\lambda$.  The developing map 
     is injective on maximal disks, and is therefore injective on the 
     homotopy class of $\gamma$.  
     
     Additionally, the extension of the developing map to the closure 
     of $\uc{S}$ is injective on the closure of each maximal disk.  
     The holonomy of $\gamma$ is hyperbolic since $S$ is closed, 
     therefore the developed image of 
     the endpoints are distinct.  Thus, $\gamma$ is 
     admissible.  \qed

 The following theorem gives a topological way to distinguish any
two complex projective structres with real holonomy. 

\begin{thm} \label{real lamination} Fix an oriented surface $S$ and a
 discrete real representation 
 \[ \rho : 
 \pi_1(S) \rightarrow \Gamma < PSL_{2}(\mathbb{R})
 .\]
 Let $\Sigma_{0} = (D_{0}, 
 \rho)$ and $\Sigma_{1} = (D_{1}, \rho)$.  Then 
 \[ \lambda_{0} \simeq \lambda_1
\Longleftrightarrow \Sigma_0 = \Sigma_1.
\]

\end{thm}

\noindent \textbf{Proof : }  
Assume $\Sigma_0 = \Sigma_1$. 
Then $D_0 \simeq D_1 $ and there is
a homeomorphism $g: S \rightarrow S$ isotopic to the
identity and a lift $\uc{g} : \uc{S} \rightarrow \uc{S}$ so
that \[ D_0 =  D_1 .\]     Now suppose $x_0 \in \uc{\lambda}_0$.
Then \[D_1  \circ \uc{g}(x_0) = D_0(x_0) \in \RP\] which
implies $\uc{g}(x_0) \in \lambda_1$ and thus \[ \uc{g}(\lambda_0)
\subset \lambda_1.\]   By switching the roles of $D_0$ and
$D_1$ we obtain an inverse, $\uc{g}^{-1}$, for $\uc{g}$.  It follows
that
$\uc{g}^{-1}(\lambda_1) \subset \lambda_0$, and consequently \[
\uc{\lambda}_1 \subset \uc{g} (\uc{\lambda}_0) \] which implies
$\uc{g} (\uc{\lambda}_0) = \uc{\lambda}_1$ and thus $g(\lambda_0) =
\lambda_1$ where $g: S \rightarrow S$ is isotopic to the identity and 
lifts to $\uc{g}$. Thus $\lambda_{1} \simeq \lambda_{2}$. The other direction of the theorem will use the 
following lemmas.  A proof of the first can be found in \cite{CB}.

\begin{lem} \label{Poi} A homeomorphism $h$ of $S$
is isotopic 
    to the identity if and only if $h$ has a lift $\uc{h}: \uc{S} 
    \rightarrow \uc{S}$ which extends to the identity on the boundary 
    of $\uc{S}$.
	

\end{lem}


\begin{lem}\label{real holonomy orientation}
    Let $\Sigma_{0} = (D_{0},\rho)$ and $ \Sigma_{1} = (D_{1}, \rho)$
    be projective structures with real discrete
    holonomy on an oriented surface $S$.  Suppose  
    that an essential nonannular subsurface $X \subset 
    S$ is disjoint 
    from $\lambda_{0}$ and $\lambda_{1}$.  Then 
    if $\uc{X}$ is a lift of $X$ in $\uc{S}$,
    $D_{0}$ and $D_{1}$ both map $\uc{X}$ into the same half plane

\end{lem}
    
\noindent \textbf{Proof : }

By Lemma \ref{max disks}
each developing map $D_i$ is injective on each component of $\uc{S} 
\backslash \uc{\lambda_i}$. 
Since $\tilde{X}$ is a contained in a component of $\tilde{S}\backslash \lambda_i$, $D_i$ is injective on $\tilde{X}$.
Let $l$ be a component of $\partial \tilde{X}$ and let $g \in \pi_1(S)$ be the deck transformation that fixes $l$. We continuously extend $D_i$ to the endpoints $p^\pm$ of $l$ in the compactification of $\tilde{S}$ by setting
$$D_i(p^\pm) = \lim_{k \rightarrow \pm \infty} \rho(g^k)D_i(x).$$
The extension of the $D_i$ to the boundary points depends only on the holonomy, not the particular developing map so the $D_i$ agree on the boundary points of $\uc{X}$.  

The orientation on $S$ lifts to an orientation on $\uc{S}$.  The orientation of $\uc{S}$ induces an orientation of the $\partial{\uc{S}}$ which induces a cyclic ordering of points in $\partial{\uc{S}}$.  In particular there is a cyclic order on the boundary points of $\tilde{X}$ in $\partial S$. Additionally, the orientation of the Riemann sphere induces an orientation on the upper and lower half planes. The orientation of the half planes then determines a cyclic ordering of the points in $\R \cup \{\infty\}$, and these two orderings are opposite.

If $D_i(\uc{X}) \subset \Omega_+$ ($\Omega_-$) then the cyclic ordering of the endpoints of the image must agree with the cyclic ordering of $\R \cup \{\infty\}$ induced by $\Omega+$ ($\Omega_-$), since $D_i$ has real holonomy and is orientation preserving.  
Since $X$ is nonannular $\uc{X}$ has at least three endpoints and thus each image $D_i(\uc{X})$ then has at least three endpoints.  Reversing the cyclic order on a set with more than three points gives a different cyclic order.
For these endpoints to have the same cyclic ordering both images must lie in the same half plane.  Therefore $D_0$ and $D_1$ must map $\uc{X}$ to the same half plane.

    \qed 


\begin{lem}\label{UL}
 Let $(D_{0},\rho)$ and $(D_{1},\rho)$ be projective structures on 
 $S$ with the same real lamination $\lambda$.  Then 
 \[ D_{0}(\uc{X}) = D_{1}(\uc{X}) \]
 where $\uc{X}$ is any component of $\uc{S} \backslash \uc{\lambda}$.

\end{lem}

\noindent \textbf{Proof : }  
By Proposition \ref{max disks} every component of $\tilde{S} \backslash \tilde\lambda$ is mapped onto either the upper or lower half plane. If $X$ is nonannular the Lemma follows from Lemma \ref{real holonomy orientation}.  Since $S$ is closed and has genus greater than 2 there is
some nonannular subsurface $Y \subset S \backslash \lambda$.
Let $X$ be a possibly nonannular component of $S \backslash \lambda$. By Lemma \ref{closed curves} there are a finite number of components $X_j \subset S \backslash \lambda$ with $j \in [0,n]$ such that $X_0 = Y$, and $X_i$ and $X_{i+1}$ are adjacent.

We now proceed by induction on the length of the sequence. We have already observed that $D_0(\tilde{Y}) = D_1(\tilde{Y})$ for all preimages $\tilde{Y}$ of $Y$. This takes care of the base case when $n=0$. Now assume that the lemma holds for all components that can be connected to $Y$ by a sequence of length $\leq n-1$ and we will prove that it holds for sequences of length $n$. If $\tilde{X}$ is a component of the preimage of $X = X_n$ then there is a component $\tilde{X}_{n-1}$ of the preimage of $X_{n-1}$ that is adjacent to $\tilde{X}$. Since $D_i$ is a local homeomorphism and $\tilde{\lambda}$ is mapped to $R \cup \{\infty\}$, adjacent components are mapped to opposite half planes. Since $D_0(\tilde{X}_{n-1}) = D_1(\tilde{X}_{n-1})$ we have $D_0(\tilde{X}) = D_1(\tilde{X})$ completing the induction step and the proof.

\qed

%
\begin{lem}\label{pushpull}
 Let $(D_{0},\rho)$ and $(D_{1},\rho)$ be projective structures on 
 $S$ with the same real curves $\lambda$.  There 
    is a $\pi_{1}(S)$-equivariant map $\uc{f} : \uc{S}
\rightarrow \uc{S}$ which is a lift of a map $f$ that is isotopic to 
    the identity such that
    \[ D_0 = 
      D_1 \circ \uc{f}. \]

\end{lem}
\noindent \textbf{Proof : }  

We first construct the map $\uc{f}$ on $\uc{\lambda}$.  This requires that
 $D_i$ for $i \in \{0,1\}$ be injective on components of $\uc{\lambda}$, and we show this now.
Let $\tilde{\gamma}$ be a component of $\tilde\lambda$ and let $\tilde{X}$ be a component of $\tilde{S} \backslash \tilde\lambda$ such that $\tilde\gamma$ is in $\partial \tilde{X}$.

  Let $\{p_1,p_2\}$ denote two points of $\tilde\gamma$ so that $D_i(p_1) = D_i(p_2)$.  Since $D_i$ is a local homeomorphism there are open neighborhoods $U_1$ and $U_2$ of $p_1$ and $p_2$ and a geometric disk $V$ in $\chat$ with center at $D_i(p_1) = D_i(p_2)$ for which $D_i$ is a homeomorphism of each $U_i$ onto $V$.  
By Lemma \ref{UL} we may assume $D_i(\uc{X}) = \Omega_+$ and we let $\{z_k\}$ be a sequence of points in $\Omega_+ \cap V$ which converge to $D_i(p_1) = D_i(p_2)$.


%

Let $x_k = D^{-1}_i(z_k) \cap U_1$ and $y_k = D^{-1}_i(z_k) \cap U_2$.  Since $D_i$ is a local homeomorphism $D^{-1}_i$ is continuous so $z_k \in \Omega_+$ implies that $x_k$ and $y_k$ are contained in $\uc{X}$.  Also, since $z_k \rightarrow D(p_1) = D(p_2)$ the continuity of $D^{-1}_i$ implies that $x_k \rightarrow p_1$ and $y_k \rightarrow p_2$.  
But Lemma \ref{max disks} implies that $\uc{X}$ is a maximal disk, so $D_i$ is injective on $\uc{X}$ and therefore $p_1 = p_2$.  This implies that $D_i$ is injective on components of $\uc{\lambda}$.  

Next we show that $D_0$ and $D_1$ map $\uc{\gamma}$ to the same arc.  As in the previous Lemma, if $p^\pm$ are the endpoints of $\uc{\gamma}$, we can extend $D_i$ to $p^\pm$ and $D_0(p^\pm) = D_1(p^\pm)$.  Since $\uc{\gamma} \subset \uc{\lambda}$ either $D_0$ and $D_1$ map $\uc{\gamma}$ to the same arc or the closure of their images is all of $\R \cup \{\infty\}$.  Let $\uc{X}$  denote a connected component of  $\uc{S} \backslash \uc{\lambda}$ whose boundary contains $\uc{\gamma}$. 
If $\overline{D_0(\uc{\gamma}) \cup D_1(\uc{\gamma})} = \R \cup \{\infty\}$, then of $C_x$ is mapped to opposite half planes by $D_0$ and $D_1$.  This is a contradiction of Lemma \ref{UL}, and so we have $D_0(\gamma) = D_1(\gamma)$.  So if $x \in \uc{\lambda}$ and $\uc{\gamma} \subset \uc{\lambda}$ contains $x$, the map \[ \uc{f}(x) = D_1^{-1}(D_0(x)) \cap \uc{\gamma} \] is well defined.

Now we extend $\uc{f}$ to all of $\uc{S}$.  For any point $x \in \uc{S} \backslash \uc{\lambda}$, let $C_x$ denote the connected component of $\uc{S} \backslash \uc{\lambda}$ containing $x$.  By Lemma \ref{max disks} the $D_i$ are homeomorphisms from $C_x$ to $\Omega_+$ or $\Omega_-$.  For any $x \in \uc{S}$ there is a unique point $ y \in
  C_{x}$ such that $D_0(x) = D_1(y)$.
 Therefore we
  define for every $x \in \uc{S}$, a map $\uc{f} : \uc{S} \rightarrow
  \uc{S} $ by \[ \uc{f} (x) = ( D_1)^{-1} (D_0(x)) \cap C_{x} .\]
    %
    %
    
    The continuity of $D_i$ implies that $\uc{f}$ is 
    continuous on each component of $\uc{S} \backslash \lambda$.  
    Since $D_i$ is a local homeomorphism, if 
    $\{x_{i}\} \in C_{x}$ and $\{y_{i}\} \in C_{y}$ converge to $x 
    \in \uc{\lambda}_{0}$ then $\{D_{0}(x_{i})\}$ and $\{D_{0}(y_{i})\}$ 
    converge to some point $r \in \R$.  Then 
    $\{D^{-1}_{1}(D_{0}(x_{i}))\} \cap C_{x}$ and 
    $\{D^{-1}_{1}(D_{0}(y_{i}))\} \cap C_{y}$ have subsequences which
    converge to $D^{-1}_{1}(r) \in \uc{\lambda}$.  Therefore $\uc{f}$ 
    is continuous on all of $\uc{S}$.  Then, for any $x \in \uc{S}$ \[ D_1 \circ \uc{f}(x) = D_1 \circ ((D_1)^{-1} (D_0(x))) \cap C_{x} = D_0(x) .\] 
    %
    %
    %
   Next we show that $\uc{f}$ is equivariant with respect to $\pi_i(S)$.   For any $g \in \pi_{1}(S)$ the equivariance of the developing map 
 implies

 \[ \begin{array}{ll}
	\uc{f}(gx) 
	& =  D_{1}^{-1}(D_{0}(gx)) \cap C_{gx} \\
 	& = D_{1}^{-1}(\rho(g)(D_{0}(x))) \cap C_{gx} \\ 
	&  = g (D_{1}^{-1}(D_{0}(x)) \cap C_{x}) = g \uc{f}(x)
	\end{array}
\] and therefore $\uc{f}$ is equivariant 
 with respect to the action of $\pi_{1}(S)$, and thus $\uc{f}$ descends to a map $f:S \rightarrow S$.  Finally, since the $D_i$ are equivalent on the endpoints of $\lambda$ the map $\uc{f}$ extends to the identity on the boundary of $\uc{S}$.  By Lemma \ref{Poi} $f$ is isotopic to the identity.

\qed



We now complete the proof of Theorem \ref{real lamination}.  Let $\Sigma_0$ and $\Sigma_1$ be projective structures with holonomy $\rho$.  By definition, $\Sigma_i$ is an equivalence class of pairs $(D_i, \rho)$ where $D_i$  is unique up to precomposition with a lift of a homeomorphism of $S$ that is isotopic to the identity and postcomposition with a Mobius transformation and $\rho$ is unique up to conjugation.  Let $D_0$ and $D_1$ be developing maps of $\Sigma_0$ and $\Sigma_1$ such that $\lambda_0$ is homotopic to $\lambda_1$.  There is a map $h_1: S \rightarrow S$ isotopic to the identity for which $h_1(\lambda_0) = \lambda_1$. Let $\tilde{h}_1$ be a lift of $h$ to $\tilde S$. Then $D'_0 = D_0 \circ \tilde{h}_1$ is equivalent to $D_0$.

 By Lemma \ref{pushpull} there is a map $h_2 : S \rightarrow S$ isotopic to the identity which has a lift $\uc{h}_2$ such that $D'_0 =  D_1 \circ \uc{h}_2$ is equivalent to  $D_1$.  Then $D_1$ and $D_0$ are in the same equivalence class and therefore $\Sigma_0 = \Sigma_1$.



\begin{notation}
   Let $\Sigma = (D,\rho)$ be a projective structure on $S$. 
   By Theorem \ref{real 
   lamination} we may denote $\Sigma$ by $\Sigma(\rho, \lambda)$ (or 
   alternatively $(\rho, \lambda)$.
   \end{notation}
 

\begin{definition}
	Let $\gamma$ be a simple closed curve in $S$.  Denote by $[\gamma]$ the homotopy class of $\gamma$.  If 
$\delta$ is an admissible curve in $[\gamma]$ then we say $\delta$ is an admissible representative of $[\gamma]$.
\end{definition}

\begin{prop} \label{disjoint grafting}
 Let $\Sigma(\rho,\lambda)$ be a projective structure with a real representation $\rho$.  If $\gamma$ is an admissible representative of $[\gamma]$ 
 disjoint from the real curves $\lambda$ of $\Sigma(\rho,\lambda)$  then  
   \[ Gr_{\gamma}(\Sigma(\rho, \lambda)) = \Sigma(\rho, \lambda \cup 2 
   \gamma)\] where $2 \gamma$ denotes two simple closed curves isotopic 
   to $\gamma$.
\end{prop}
 
 \noindent \textbf{Proof : }  
 By definition, grafting $\Sigma(\rho,\lambda)$ along $\gamma$ produces a new
projective structure $(D_1,\rho) = \Sigma(\rho, \lambda_{\gamma})$.
Since grafting affects the structure locally near $\gamma$, the fact that 
$\gamma$ is disjoint from $\lambda$ implies that $D_{1}(\uc{\lambda}) = D(\uc{\lambda}) = \RP$ and so 
$\lambda \subset \lambda_{\gamma}$.

Let $A$ denote the annulus grafted into $\Sigma(\rho,\lambda)$ along $\gamma$, and let $L$ be a component of $\uc{A}$.  There is some $g \in \pi_1(S)$ which fixes $L$ and thus $\rho(g)$ fixes the annulus $D_{1}(L)$  in $\chat$.  The fixed points $\{p_1,p_2\}$ of $\rho(g)$ form the boundary of $D_1(L)$ and split 
$\R \cup \{\infty\}$ into two components that are both contained in $D_{1}(L)$.    
Then, in $L$ there 
are two arcs $\uc{\gamma}_{1}$, $\uc{\gamma}_{2}$
of $D_{1}^{-1}(\RP)$ and the boundaries of these arcs are the fixed points of $g$.  Therefore
$\gamma_i = \uc{\gamma}_{i} \slash g$ are simple 
closed curves in $A$ isotopic to $\gamma$ and \[ \lambda_{\gamma} 
= \lambda \cup \gamma_{1} \cup \gamma_{2} .\] \qed

Assume $\alpha$ and $\beta$ are admissible representatives for $[\alpha]$ and $[\beta]$ in a projective structure
$\Sigma = (\rho, \lambda)$, and \[\alpha \cap \beta = \alpha \cap \lambda = \beta \cap \lambda = \phi.\]  By Proposition \ref{disjoint grafting}, $Gr_\alpha (\rho,\lambda) = (\rho, \lambda \cup 2 \alpha)$ and $Gr_\beta (\rho, \lambda) = (\rho, \lambda \cup 2\beta)$.  By Lemma \ref{disjoint admissibility} $\alpha$ has an admissible representative in $(\rho, \lambda \cup 2 \beta)$ and $\beta$ has an admissible representative in $(\rho, \lambda \cup 2 \alpha)$.  Proposition \ref{disjoint grafting} again implies $Gr_{\alpha}(\rho, \lambda \cup 2 \beta) = (\rho, \lambda \cup 2 \beta \cup 2 \alpha)$ and $Gr_{\alpha}(\rho, \lambda \cup 2 \alpha) = (\rho, \lambda \cup 2 \beta \cup 2 \alpha)$.  Then Theorem
\ref{real lamination} implies \[Gr_{\alpha}(Gr_\beta(\rho,\lambda)) = 
Gr_{\beta}(Gr_\alpha(\rho,\lambda)).\]  

Let $\sigma = \{\alpha_{i}\}$
denote a (possibly empty) collection of disjoint admissible curves in $(\rho, \lambda)$ which are each disjoint from $\lambda$.  By the above argument we may denote \[Gr_{\sigma}(\Sigma) = Gr_{\alpha_{1}} \circ
Gr_{\alpha_{2}}\circ \dots \circ Gr_{\alpha_{n}}(\Sigma))\]

\section{Fuchsian Holonomy}

%
%
In this section we use ideas from the previous section to prove a theorem of Goldman which classifies projective structures with Fuchsian holonomy \cite{Gol}.  Our proof is similar to Kapovich's proof in \cite{Kap89} of a generalization of this theorem to all dimensions.

\begin{definition}
       A projective structure $\Sigma$ is standard if its developing map
is a covering map onto its image.  Equivalently, the projective
structure is the quotient of the domain of discontinuity by the
holonomy representation \[ \Sigma = \Omega \hs \slash \rho. \] 
\end{definition}

\begin{thm}[Goldman]\label{Goldman}
If $\Sigma = (\rho, \lambda)$ is a complex projective structure
with $\rho$ a Fuchsian representation, then there exists a collection 
of admissible disjoint simple closed curves $\sigma$ such that \[ \Sigma = Gr_{\sigma}(\Sigma_{0}) \] where 
$\Sigma_{0}$ is a standard (hyperbolic) complex projective structure.
    \end{thm}
    
    \noindent \textbf{Proof : }  
	Any Fuchsian projective structure $\Sigma = (\rho, 
	\lambda)$ has the property that each simple closed curve
	in $\lambda$ contains an even number of curves in its isotopy 
	class.  This follows from Lemma \ref{real holonomy orientation} in the following 
	way. 

 Let $\Sigma_{0} = (\rho, \lambda)$ be a standard Fuchsian projective 
	structure.  Let $U \subset S$ be any non-annulur region 
	disjoint from $\lambda$.  Since $\Sigma_{0}$ is standard and $\rho$ is Fuchsian, $U$ 
	is disjoint from $\lambda_{0} = \phi$. Now Lemma \ref{real 
	holonomy orientation} implies that $D(U)$ and $D_{0}(U)$ are both 
	positive (contained in $\Omega_{+}$) or negative  
	(contained in $\Omega_{-}$).
	
	This imples that every non-annulur region in $S$ that is 
	disjoint from $\lambda$ is, say, positive.  Since adjacent regions of $S \backslash \lambda$ are mapped to opposite half planes the number of
	disjoint, simple closed isotopic curves in $\lambda$ must be 
	even.  Since $S$ is compact the number of such curves in $\lambda$ is finite.
	
	Enumerate the real curves of $\Sigma$ as $\lambda = 2 \gamma_1 \cup 2 \gamma_2 \cup \dots \cup 2 \gamma_n$.  Let 
	$\sigma = \{\delta_{1} \dots \delta_{n}\} $ denote a collection of $n$ disjoint curves in $S$, 		where $\delta_{i}$ is isotopic to $\gamma_{i}$.  Since $\Sigma_0$ is standard $\delta_i \cap \lambda_0 = \phi$ so Lemma \ref{disjoint admissibility} implies that $\delta_i$ is admissible.  
	By Proposition \ref{disjoint grafting} \[ Gr_\sigma (\Sigma_0) = (\rho, 2 \delta_1 \cup 2 \delta_2 \cup  \dots \cup 2 \delta_n) \] and Theorem \ref{real lamination} then implies $Gr_\sigma(\Sigma_0) = \Sigma$. \qed
	

\section{Real Schottky Holonomy} \label{Real Schottky Holonomy}

In this section we consider complex projective structures whose 
holonomy representation is a Schottky subgroup of 
$PSL_{2}(\R)$.  Since these representations are quasi-conformally conjugate to real representations, we will use ideas developed in the previous section.  We define a Schottky group as in \cite{Kap}.

\begin{definition} \label{Schottky def}
     Let $B_1, \dots, B_k, B'_1, \dots, B'_k$ be disjoint closed geometric disks in $\chat$, and let int($B_i$) denote the interior of $B_i$ and ext($B_i$) denote the exterior of $B_i$.   Let $g_i$ be a    
     Mobius transformation such that $g_i$: int($B_i$) $\rightarrow$ ext($B'_i$) .
    The group $\Gamma = < g_1, g_2, \dots, g_k>$ is called a classical Schottky 
     group.  
\end{definition}

Since $\rho$ is quasi-conformally conjugate to a representation from $\pi_1(S) \rightarrow \Gamma \subset PSL_2(\R)$ it suffices to prove the case where $\rho$ is a real Schottky representation.  Unless otherwise stated, from here forward  we fix a representation $\rho: \pi_{1}(S) \rightarrow \Gamma$ where $\Gamma < PSL_{2}(\mathbb{R})$ is  a real Schottky group.  We 
   also assume that all curves intersect transversely.  
   Theorem \ref{real 
   lamination} implies that our study of grafted real
   projective structures is 
   reduced to a study of the behavior of their real curves near 
   the grafting curve.  With this in mind we begin with 
   a careful examination of the universal cover of the Hopf torus and its intersection with 
   the real line.

\section{Basic Facts for Curves in a Torus}

Let $(\alpha,\beta)$ be an ordered pair of essential oriented transverse curves in a torus $T$ which intersect exactly once.  Then $(\alpha,\beta)$ determines an orientation for $T$ and we choose this orientation for $T$. 
Given any orientation of $T$ we define the index of an intersection point of two curves $\alpha$ and $\beta$ to be $+1$ if the orientation on $T_x(T)$ induced by $(\alpha, \beta)$ agrees with the given orientation of $T$, and $-1$ otherwise.  The algebraic intersection number $\hat{i}(\alpha,\beta)$ is defined to be the sum of the indices of all the intersections of $\alpha$ and $\beta$.  
Since $H_1(T) \simeq \mathbb{Z} \oplus \mathbb{Z}$, given a basis for $H_1$ the homology class of any simple closed multicurve is given by an ordered pair of integers $(p,q)$.  An orientation on a simple closed curve in $T$ is a choice of sign of the associated integer pair thus $(p,q)$ has the opposite orientation of $(-p,-q)$.  

\section{The Hopf Torus} \label{section: the hopf torus}
\begin{figure}[ht]
	\begin{center}
	  \psfrag{1}{$\mathbb{R} \subset \chat$}
	  \psfrag{2}{$\hat{\gamma}$}
 	 	\includegraphics[scale=.7]{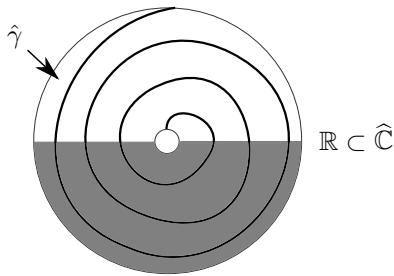}
	    \caption{ The developed image of $\gamma$ when $\gamma$ is left-spiraling.}
	    \label{left-spiraling}
	\end{center}
	    \end{figure}

Next we consider a particular class of admissible curves in $S$ whose associated graftings can be described by the image of these curves in the Hopf torus.  A lift of such a curve spirals around its fixed points as suggested in Figure \ref{left-spiraling} and we call these curves spiraling curves.  This reduces the study of projective structures grafted along such curves to a study of curves in a torus.

%
%

Let $\gamma_T = \hat{\gamma} \slash \rho(\gamma)$ denote the projection of $\gamma$ to the Hopf torus.  We claim that $\gamma_T$ has minimal intersection with $\RP \slash \rho(\gamma)$.   Assume $\gamma_T$ does not have minimal intersection with $\RP \slash \rho(\gamma)$.  Then there is an innermost disk in $T_H$ bounded by $\RP \slash \rho(\gamma)$.  Since $\R^-$ and $\R^+$ are disjoint,  any such disk lifts to a disk in $\chat$ whose boundary is composed of two arcs, one from $\hat{\gamma}$ and one from either $\R^-$ or $\R^+$.  Then consecutive points of $\hat{\gamma} \cap \RP$ lie on a single component of $\RP \backslash \{0,\infty\}$ and $\gamma$ cannot be spiraling.   

 Recall, the Hopf torus $T_H$ is defined as the quotient  $(\chat \backslash Fix(\rho(\gamma)) \slash \rho(\gamma)$.  So any curve in $S$ with holonomy $\rho(\gamma)$ projects to an essential closed curve in $T_H$.  Let $\alpha$ be a geometric circle in $\chat$ which separates $Fix(\rho(\gamma))$.   Since $Fix(\rho(\gamma)) = \{0,\infty\}$ it follows that $| \alpha \cap \R^-| = 1$ and thus $|\alpha_T \cap \lambda_-| = 1$ where  $\alpha_T = \alpha \slash \rho(\gamma)$ and $\lambda_- = \R^- \slash \rho(\gamma)$.  

We orient $\alpha_T$ and $\lambda_-$ so that 
that the orientation of $T$ produced by the oriented ordered pair $(\lambda_-,\alpha_T)$ agrees with the standard orientation of $T$.  We choose as a basis for $H_1(T_H)$ the ordered pair of oriented simple closed curves $(\lambda_-,\alpha_T)$.  Now $\lambda_{-} = (1,0)$ and $\alpha_T = (0,1)$. 

Since $\alpha$ separates the endpoints of $\hat{\gamma}$ we have $\hat{i}(\hat{\gamma},\alpha) = 1$.  So in the Hopf torus $\hat{i}(\gamma_T,\alpha_T) = 1$. Since $\gamma_T$ is closed \[ \gamma_T = (1,k) \]  for $k \in \mathbb{Z}$ and we say $\gamma_T$ is left-spiraling if $k > 0$ and is right-spiraling if $k < 0$.

\section{Flat and Sharp Operations} 

\begin{figure}[ht]
     \psfrag{1}{\small{$\sharp$}} \psfrag{2}{\small{$\flat$}}
	    \centerline{ \epsfig{file=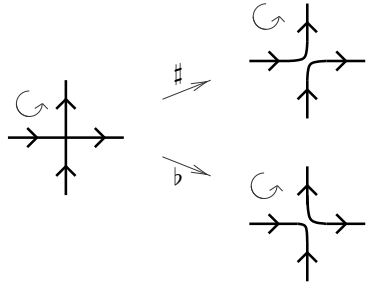,scale=.8} }
	       \caption{Flat and sharp operations}
	       \label{sharp flat fig}
  \end{figure}

Next we define a surgery operation on curves in a surface.  We will 
later express the real curves obtained by grafting along spiraling
curves as the result of this surgery on  
the real curves of the initial structure together with curves 
isotopic to the grafting curve.  The surgery is depicted in Figure \ref{sharp flat fig}.

Let $(\lambda, \gamma)$ be an ordered pair of
simple cloed multi-curves in an oriented surface $S$. For each
$x \in \lambda \cap \gamma$ let $N(x)$ denote a neighborhood of 
$x$ in $S$ so that $N(x) \cap (\lambda \cap \gamma) = x$.  
Choose a local orientation of $\lambda$ and $\gamma$ near each $x \in \lambda \cap \gamma$ 
such that the orientation induced by $(\lambda, \gamma)$ agrees with the orientation of $S$.

The multicurve $[\lambda, \gamma]_{\sharp}$ is obtained by
resolving each $x \in \lambda \cap \gamma$ by replacing $x$ with
arcs in $N(x)$ joining $\lambda$ to $\gamma$ so that the local
orientations of the arcs being joined agree.  The $\flat$
operation joins arcs whose local orientations disagree.   
If $\lambda$ and $\gamma$ are given the opposite orientation 
the same topological curves are obtained, but with the opposite orientation.

 \subsection{Operations on the Torus} \label{torus}
  

  It will often be the case that these operations on the surface $S$ may be reduced to the same operations on the Hopf torus. 
Here we derive relations between $\sharp \slash \flat$ operations in the case of the torus.
Recall that we assume all intersections are transverse and minimal.

    \begin{lem} \label{flat sharp switch} 
	    Let ($\alpha, \beta$) be an ordered pair of essential simple closed curves in an oriented surface $S$.  Then 
    	\[ [\alpha, \beta]_\sharp =  [\beta, \alpha]_\flat.\]    
    \end{lem}
        
        \noindent \textbf{Proof :}

The local orientations required to compute $[\beta, \alpha]_\flat$ may be obtained from those required to compute $[\alpha, \beta]_\sharp$ by reversing the local orientations of exactly one of the curves near each $x_i$.  
Resolving intersections so that the orientations of the curves being joined disagree is the same as first reversing the orientations of one of the curves and then resolving intersections so that orientations of the curves being joined agree.  Therefore the first and second computation produce the same set of curves. \qed

Assume $\alpha$ and $\beta$ are oriented essential simple closed curves in $T$ in minimal position. 
The $\sharp$ operation depends only on the orientation of the surface, not on the orientations of the curves themselves, but we will use the oriention $\alpha$ to define an orientation of $[\alpha,\beta]_\sharp$ and $[\alpha,\beta]_\flat$. 

Since $\alpha$ and $\beta$ are curves in a torus, each index of an intersection point of $\alpha$ with $\beta$ is the same for each point.  Choose an orientation of $\beta$ so that each index is +1.  Then by definition the $\sharp$ operation joins $\alpha$ and $\beta$ so that the orientations agree, and we choose this orientation for  $[\alpha,\beta]_\sharp$.  
Then $[\alpha, \beta]_\flat$  inherits an orientation in a similar manner if we choose an orientation of $\beta$ so that each index is -1. 

The following lemma exhibits the relationhip between algebraic intersection number and the $\sharp$ and $\flat$ operations on curves in a torus.


\begin{lem} \label{sharp respects a.i.}
    Assume that  $\alpha, \beta$  and  $\gamma$ are oriented simple closed curves in an oriented torus $T$.  Then
$\hat{i}(\alpha,\beta) > 0$ implies
    
    \[  \begin{array}{ll} 
     \hat{i}([\alpha,\beta]_{\sharp}, \gamma) = 
    	\hat{i}(\alpha, \gamma) + \hat{i}(\beta, \gamma) \\
     \hat{i}([\alpha,\beta]_{\flat}, \gamma) = 
    	\hat{i}(\alpha, \gamma) - \hat{i}(\beta, \gamma)
    \end{array}
    \]
and $\hat{i}(\alpha,\beta) < 0$ implies
    
    \[  \begin{array}{ll} 
     \hat{i}([\alpha,\beta]_{\flat}, \gamma) = 
    	\hat{i}(\alpha, \gamma) + \hat{i}(\beta, \gamma) \\
     \hat{i}([\alpha,\beta]_{\sharp}, \gamma) = 
    	\hat{i}(\alpha, \gamma) - \hat{i}(\beta, \gamma)
    \end{array}
    \]

    \end{lem}
\noindent \textbf{Proof : } 
Transversality implies that the intersection of $\alpha$ and $\beta$ is disjoint from $\gamma$.  We may then assume that each point of $\gamma \cap [\alpha,\beta]_\sharp$ corresponds exactly with either a point of $\gamma \cap \alpha$ or $\gamma \cap \beta$.  
 The assumption that $\hat{i}(\alpha,\beta) > 0$ implies that at each $p \in \alpha \cap \beta$ the orientations of the curves agree with the orientation of $T$.  Therefore the orientation of $[\alpha,\beta]_\sharp$ is the same as $\alpha$ where the they coincide, and similarly for $\beta$.  Therefore \[ \hat{i}([\alpha,\beta]_{\sharp}, \gamma) = 
\hat{i}(\alpha, \gamma) + \hat{i}(\beta, \gamma) .\]  

 Since $[\alpha, \beta]_\sharp = [\alpha, -\beta]_\flat$ it follows from $\hat{i}(- \beta, \gamma) = -\hat{i}(\beta, \gamma)$ that  \[ \hat{i}([\alpha,\beta]_{\flat}, \gamma) = 
\hat{i}(\alpha, \gamma) - \hat{i}(\beta, \gamma) .\] 

If $\hat{i}(\alpha,\beta) < 0$ then $\hat{i}(\alpha,-\beta) > 0$ so the argument above implies   
\[ \hat{i}([\alpha,\beta]_\flat, \gamma) =  \hat{i}(\alpha, \gamma) - \hat{i}(-\beta, \gamma) =  \hat{i}(\alpha, \gamma) + \hat{i}(\beta, \gamma) \]  

and  \[ \hat{i}([\alpha, \beta]_{\sharp}, \gamma) = \hat{i}(\alpha, \gamma) + \hat{i}(-\beta, \gamma) = \hat{i}(\alpha, \gamma) - \hat{i}(\beta, \gamma)  .\] 

\qed


In the discusssion above the ordered pair of integers $(p,q)$ denotes 
   a homology class of an oriented simple closed multicurve in a torus.  We will also use this same notation to denote a representative curve from the homology class. 

		\begin{lem} \label{arcs in torus} Let $\alpha = (p,q)$ and $\beta = (r,s)$ denote two oriented essential simple closed curves in an oriented torus $T$.  Then $\hat{i}(\alpha,\beta) \geq 0$ implies 
    \[ [(p,q),(r,s)]_{\sharp} = (p+r,q+s) \] and \[   
    [(p,q),(r,s)]_{\flat} = (p-r,q-s) \]
and $\hat{i}(\alpha,\beta) < 0$ implies
    \[ [(p,q),(r,s)]_{\flat} = (p+r,q+s) \] and \[   
    [(p,q),(r,s)]_{\sharp} = (p-r,q-s). \]

\end{lem}
    \noindent \textbf{Proof : } 
%
If $ \hat{i}((p,q),(r,s) = 0$ then the curves are disjoint and are thus homotopic, so the Lemma holds in this case.
Assume $\hat{i}((p,q),(r,s) > 0$.

By Lemma \ref{sharp respects a.i.}
\[ \hat{i}([(p,q),(r,s)]_{\sharp}, (1,0)) = 
\hat{i}((p,q), (1,0)) + \hat{i}((r,s), (1,0)) = -q - s\]  and
\[ \hat{i}([(p,q),(r,s)]_{\sharp}, (0,1)) = 
\hat{i}((p,q), (0,1)) + \hat{i}((r,s), (0,1)) = p + r .\] 	
	
Then 
$[(p,q),(r,s)]_{\sharp}  = (p+r, q+s)$. Since $[\alpha, \beta]_{\flat} = [\alpha, -\beta]_{\sharp}$ it follows that
\[[(p,q),(r,s)]_{\flat} = ((p,q),-(r,s)]_{\sharp} = 
[(p,q),(-r,-s)]_{\sharp} = (p-r,q-s) .\]  A similar proof works for the case that $\hat{i}((p,q),(r,s) < 0$.  \qed

    
%
%
%
%


\section{Proof of the Main Theorems} \label{main proofs}

\begin{figure}[ht]
	\begin{center}
	  \psfrag{1}{$\lambda$}
	  \psfrag{2}{$\lambda_\gamma = [\lambda,2 \gamma]_\flat$}
	  \psfrag{3}{$\lambda_{\gamma_T} = [\lambda_T,2 \gamma_T]_\flat$}
	  \psfrag{4}{$\lambda_T$}
	  \psfrag{5}{$2 \gamma_T$}
	  \psfrag{6}{$A = N(\gamma)$}
	  \psfrag{7}{$\psi$}

 \includegraphics[scale=.6]{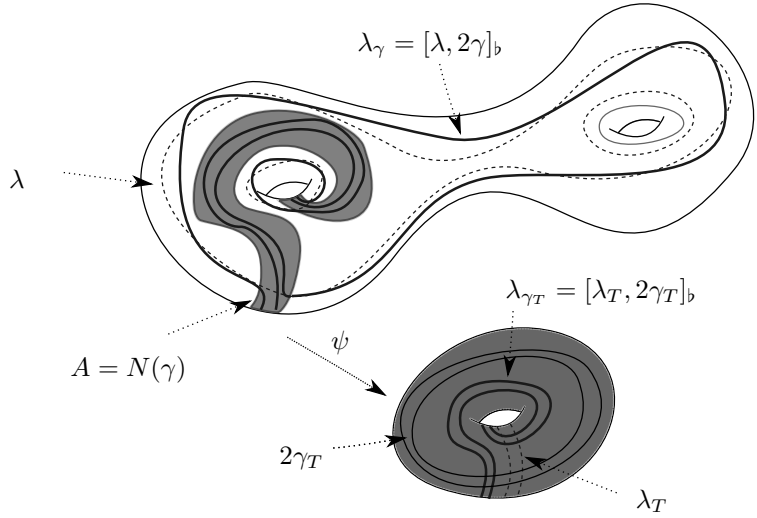}
 
  \caption{The quotient of the grafting annulus by its boundary.}
  \label{fig: quotient}
	\end{center}
	\end{figure}

Figure \ref{fig: quotient} illustrates the setting of the following technical lemma, which will be used to prove the main theorem.

\begin{lem} \label{sharp passes}
    Let $\lambda$, $\lambda_\gamma$ and $\gamma$ be oriented simple closed curves in an oriented 
    surface $S$.  Assume $\lambda_\gamma$ and $\lambda$ are homotopic relative their boundary in the complement  		of an 
    annulus $A$ which is homeomorhpic to a regular neighborhood of 
    $\gamma$.  Let $\psi: A \rightarrow T$ be a map onto a torus
    which identifies the boundary 
    components of $A$ so that the image under $\psi$ of each component of $\lambda_\gamma \cap A , \lambda \cap A$ and $\gamma \cap A$ is a simple closed curve in $T$.
    Then \[ \lambda_{\gamma_{T}} \simeq [\lambda_{T}, 2 \gamma_{T}]_{\sharp}
    \Leftrightarrow \lambda_\gamma    \simeq [\lambda, 2 \gamma]_{\sharp} \] and     
    \[ \lambda_{\gamma_{T}} \simeq [\lambda_{T}, 2 \gamma_{T}]_{\flat}
       \Leftrightarrow \lambda_\gamma \simeq [\lambda, 2 \gamma]_{\flat} \]

	where $\lambda_{\gamma_T} = \psi (\lambda_\gamma \cap A), \lambda_T = \psi (\lambda \cap A)$ and $\gamma_T = \psi (\gamma \cap A)$.
    
    \end{lem}
    
    \noindent \textbf{Proof : } 
	By symmetry of the $\flat \slash \sharp$ operations we may work with just one, say $\flat$. 
Let $2 \gamma$ denote two disjoint curves isotopic to $\gamma$, and 
assume $2 \gamma \subset A$, which implies all intersections of
$2 \gamma \cap \lambda$ 
are contained in the interior of $A$.  

The $\flat$ operation affects $2 \gamma$ and $\lambda$ only near their intersection and $\psi$ is a quotient map on $A$ which defines $T$ so it follows that  
\[ ([\lambda, 2 \gamma]_{\flat} \cap A) \slash \psi =
[\lambda_{T},2 \gamma_{T}]_{\flat}.\]  

First we assume $[\lambda_{T}, 2 
\gamma_{T}]_{\flat} \simeq \lambda_{\gamma_{T}} $.  Then by definition of $\lambda_{\gamma_T}$,  $[\lambda_T, 2 \gamma_T]_\flat \simeq (\lambda_\gamma \cap A) \slash \psi$, and the statement above implies \[([\lambda, 2 
\gamma]_{\flat} \cap A) \simeq \lambda_\gamma \cap A\] where homotopy here is relative to the boundary.  And since $\lambda_\gamma \simeq \lambda$ outside of $A$ this now implies \[ [\lambda, 2 \gamma]_\flat \simeq \lambda_\gamma.\] 
Now assume $ \lambda_\gamma    \simeq  [\lambda, 2 \gamma]_{\flat}$.  The assumption that $\lambda_\gamma \backslash A  \simeq \lambda \backslash A$  now implies that $\lambda_\gamma \cap A \simeq [\lambda, 2 
\gamma]_{\flat} \cap A$ where homotopy is relative to the boundary.  Therefore $\lambda_{\gamma_{T}} \simeq [\lambda_{T}, 2 
\gamma_{T}]_{\flat}$. Then \[ (\lambda_\gamma \cap A) \slash \psi \simeq ([\lambda,2 \gamma]_\flat \cap A) \slash \psi = [\lambda_T,2 \gamma_T]_\flat.\]

\qed

%
%

\subsection{The Hopf Torus}

To prove our main theorem, below, we use Theorem \ref{real lamination} to classify (spiraling) graftings of projective structures in terms of essential simple closed curves of the original projective structure.   Our strategy is  to calculate how the real curves of the original structure change under grafting along spiraling curves. Lemma \ref{sharp passes} shows that it suffices to do the calculation inside the grafting annulus.  To make this calculation we need the following definition, which allows us to speak of a curve before and after grafting.  
 
\begin{definition} \label{preimage}
 
Let $\delta$ be a closed curve in  $\Sigma(\rho, \lambda)$.  The expansion of $\delta$ under grafting along a curve $\gamma$, is defined as 
 \[ \delta_{\gamma} = \nu^{-1}(\delta) \] where $\nu$ is the collapsing 
 map defined in Section \ref{collapse}.
	 
\end{definition}

There is a homotopy equivalence of $S$ which identifies every curve in $S$ with its expansion under grafting.  However, the curves are geometrically different inside the grafting annulus and a precise calculation of this difference when the grafting curve is spiraling is used to prove our main result, Theorem \ref{flatsharp}.  In the next paragraph we develop notation which we will use in the proof of the main theorem.

\subsubsection{Real Curves in the Hopf Torus} \label{subsubsection: RCitHT}

\begin{figure}[ht]
	\begin{center}
	  \psfrag{2}{{$S$}}
	  \psfrag{1}{{$\lambda' \cap A$}}
	  \psfrag{3}{$\lambda_\gamma$}
	  \psfrag{0}{{$\lambda$}}

	  \psfrag{7}{$\chat$}
	  \psfrag{6}{{$\hat{\R}$}}
	  \psfrag{4}{{$D(\uc{\lambda'}) \cap \uc{A}$}}
	  \psfrag{5}{{$D(\uc{\lambda}_\gamma)$}}
	  \psfrag{12}{{$D(\uc{\gamma})$}}
	  
	  \psfrag{11}{{$\lambda_T$}}
	  \psfrag{8}{$\lambda_{\gamma_T}$}
	  \psfrag{10}{{$\gamma_T$}}
	  \psfrag{9}{$T$}

 \includegraphics[width=300px,height=200px]{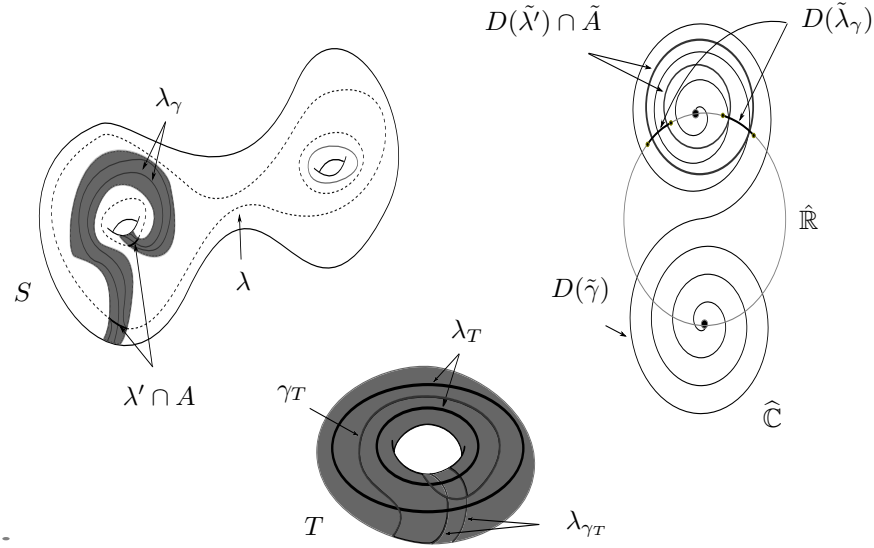}
 
  \caption{Grafting $S$ along a left-spiraling curve $\gamma$, the developed image and the Hopf torus.  }
  \label{fig: spiraling_hopf_torus}
	\end{center}
	\end{figure}

Assume $\gamma$ is an admissible
 representative in  $\Sigma(\rho, \lambda) = (D,\rho)$.  
  Let $A$ denote the grafted annulus in  $Gr_{\gamma}(\Sigma(\rho, \lambda)) = (D', \rho)$ and let $\lambda' = \nu^{-1}(\lambda)$.   Let $\uc{\lambda}'$ and $\uc{A}$ denote components of the preimage of $\lambda'$ and $A$ in $\uc{S}$ which intersect.   Let $g \in \pi_1(S)$ be the element that fixes $\uc{A}$. Since $\lambda$ is separating and $\gamma$ intersects $\lambda$ transversely set \[ \nu(\lambda' \cap A) = \{x_1, \dots, x_{2k}\} \] where $k$ is a positive integer. Recall that in the construction of $h$ in Section $\ref{collapse}$, arcs of $A$ that are collapsed to points of $\gamma$ develop to simple closed curves in $\chat \slash Fix(\rho(g))$.  This implies $\lambda' \cap A$  develops to an even number of simple closed curves in $\chat \slash Fix(\rho(g))$ and these arcs are contained in a fundamental domain of $\chat \slash Fix(\rho(g))$ for the action of $\rho(g)$.  The quotient of these curves is a multi-curve in the Hopf torus  \[ \lambda_T = D'(\uc{\lambda}' \cap \uc{A}) \slash \rho(g)  = (0, 2k) \label{concentric}  \] 
with $k > 0$ in a basis $(\lambda_-,\alpha_T)$ for $T$ with $\lambda_-$ and $\alpha_T$ as defined in Section \ref{section: the hopf torus}. Figure \ref{fig: spiraling_hopf_torus} describes the situation where $\gamma = (1,1)$ is a left spiraling admissible curve and the real curves of the grafted structure are denoted by $\lambda_\gamma$.  This figure will help to explain the proof of Theorem \ref{flatsharp}.

 \begin{thm} \label{flatsharp} Let $\Sigma = \Sigma(\rho, \lambda)$ be
 a Schottky projective structure on $S$.  Let $\gamma$ be an
 admissible right-spiraling representative in $S$.  Then \[
 Gr_{\gamma}(\Sigma(\rho,\lambda)) = \Sigma(\rho,\lambda_{\gamma}) \]
 with \[ \lambda_{\gamma} \simeq [\lambda, 2 \gamma]_{\flat} \] 
 where  $2\gamma$ denotes two parallel copies of $\gamma$.
	   \end{thm}

\noindent \textbf{Proof : } 

The first step in the proof is to examine the Hopf torus and show \begin{equation} \label{eq_1} \lambda_{\gamma_T} \simeq [\lambda_T,2 \gamma_T]_\flat .\end{equation}  Our strategy for this is to calculate the left and right sides of this equation independently, and show that both equal the same homotopy class of curves on the torus.  We will then use Lemma \ref{sharp passes} to show that this equality for curves in the torus implies the equality for curves in $S$.

Recall, that since $\rho$ is quasi-conformally conjugate to a representation from $\pi_1(S) \rightarrow \Gamma \subset PSL_2(\R)$ it suffices to prove the case where $\rho$ is a real Schottky representation.  With $\lambda_-$ and $\alpha_T$ defined above, we fix $(\lambda_-,\alpha_T)$ as a basis for the Hopf torus $T$.  Also, and let $\uc{A}$ and $g$ be defined as in the preceding paragraph.

First we examine the right side of Equation (\ref{eq_1}). As discussed in Section \ref{subsubsection: RCitHT}, the quotient by $\rho(g)$ of the developed image of the expansion of $\lambda$ by grafting along $\gamma$ \[\lambda_{T} = D'(\uc{\lambda}' \cap \uc{A}) \slash \rho(g) = (0, 2k) \mbox{ , } k > 0 \] is a multicurve with an even number of components in the Hopf torus $T$, each homotopic to $\alpha_T$.   
By definition of right-spiraling $\gamma_T = (1,-k)$ in $T$ and we let $2 \gamma_{T} = (2,-2k)$ denote two parallel copies of $\gamma_T$ in $T$. 
Then since the algebraic intersection number of $\lambda_T$ and $\gamma_T$  \[\hat{i}((0,2k),(2,-2k)) = -4k\] is negative , Lemma \ref{arcs in torus} implies \[  [ \lambda_{T}, 2 \gamma_{T}]_{\flat} = (0,2k) + (2,-2k) = (2,0)  \] in $T$. %

Now we turn our attention to the left side of Equation (\ref{eq_1}).  By definition, $\lambda_\gamma$ develops to $\R$ under the new grafted developing map.  Since the holonomy is real, the fixed points of $\rho(\pi_1(S))$ are contained in $\R \cup \{\infty\}$, thus the image of $\R$ in the Hopf torus, $\lambda_{\gamma_T}$ is the multi-curve $(2,0)$ given our choice of basis (see Section \ref{section: the hopf torus}).  Therefore \[ \lambda_{\gamma_T} = h(\lambda_\gamma \cap A) = (2,0) = [\lambda_{T}, 2 \gamma_{T}]_{\flat}  \] in $T$ (see Figure \ref{fig: spiraling_hopf_torus}). We have now established Equation (\ref{eq_1}).

Next we use Lemma \ref{sharp passes} and Equation (\ref{eq_1}) to finish the proof.  Grafting affects the real curves only in a neighborhood of the 
grafting curve $\gamma$,  so $\lambda_\gamma$ and $\nu^{-1}(\lambda)$ are homotopic in the complement of $A$ (see Figure \ref{fig: quotient}).  Given  Equation (\ref{eq_1}) for curves on the Hopf torus, we now apply  
Lemma \ref{sharp 
passes} to get equality for curves on the surface \[ \lambda_\gamma = [\nu^{-1}(\lambda), 2 
\gamma]_{\flat} \simeq [\lambda, 2 \gamma]_\flat\] and the theorem is proved. 
Similarly, if $\gamma$ is left-spiraling then $2 
\gamma_{T} = (2, 2k)$,  and thus 
\[ \lambda_{\gamma} = [\lambda, 2 \gamma]_{\sharp} .\]  \qed

\begin{figure}[ht]
  \psfrag{6+}{\small{$(y^{+}_{1},z^{+}_{1})$}}
  \psfrag{6-}{\small{$(y^{-}_{1},z^{-}_{1})$}}
  \psfrag{5+}{\small{$(y^{+}_{6},z^{+}_{2})$}}
  \psfrag{5-}{\small{$(y^{-}_{6},z^{-}_{2})$}}
  \psfrag{4+}{\small{$(y^{+}_{2},z^{+}_{3})$}}
  \psfrag{4-}{\small{$(y^{-}_{2},z^{-}_{3})$}}
  \psfrag{3+}{\small{$(y^{+}_{5},z^{+}_{4})$}}
  \psfrag{3-}{\small{$(y^{-}_{5},z^{-}_{4})$}}
  \psfrag{2+}{\small{$(y^{+}_{3},z^{+}_{5})$}}
  \psfrag{2-}{\small{$(y^{-}_{3},z^{-}_{5})$}}
  \psfrag{1+}{\small{$(y^{+}_{4},z^{+}_{6})$}}
  \psfrag{1-}{\small{$(y^{-}_{4},z^{-}_{6})$}}
  \centerline{ \epsfig{file=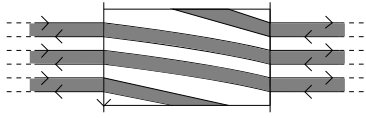,scale = .5}
  }
  \caption{A lift of the grafting annulus and preimage of the real curves}
  \label{twostep}
  \end{figure}

    
    In Figure \ref{twostep} a fundamental domain for a component of the universal cover of the grafted annulus is depicted along with its intersection with lifts of the newly obtained real curves resulting from grafting.  The horizontal curves entering the strip are lifts of the original real curves while the curves inside are lifts of the  new real curves produced by grafting.
 
 \section{Grafting Complex of Schottky Projective Structures}
 \label{grafing complex}

\begin{definition}
    We say two projective structures differ by an elementary twist if 
    their real laminations differ by a Dehn twist along a meridian curve
    which intersects one of the laminations exactly twice.
    \end{definition}
    
Let $\mathcal{G}^*(S)$ be the graph of marked Schottky projective
structures where each vertex is a structure, and an edge between two
vertices exists if the two corresponding structures differ by an
elementary twist or by grafting.  In the remainder of this section unless explicitly stated otherwise $\beta, \gamma$ and  $\lambda$ are  essential simple closed multi-curves in $S$ such that 
$| \beta \cap \gamma | = 1$, $| \beta
\cap \lambda | = 2$.  For an integer $k$ let $T^{k}_{\beta}$ denote $|k|$ Dehn twists about $\beta$, where the twists are positive or negative according to the sign of $k$.

In \cite{FM} Chapter 3 page 64, the authors observe that the effect on curves of a (right) Dehn twist is realized via a surgery operation that is equivalent to our $\flat$ operation.   If $a$ and $b$ are simple closed curves, to realize $T_a(b)$, the set of curves $i(a,b) \times a \cup  b$ is surgered so that if one follows an arc of $b$ towards the intersection, the surgered arc turns right	 at the intersection.  This is precisely our $\flat$ operation so the following relationship holds between Dehn twists and $\flat$ operation.  

\begin{lem} \label{twists and sharp} If $\alpha$ and $\beta$ are essential simple closed curves with minimal intersection in a surface $S$ then 
\[ T_\alpha (\beta) = [\beta, i(\alpha,\beta) \alpha]_\flat \] 
and
\[ T^{-1}_\alpha (\beta) = [\beta, i(\alpha,\beta) \alpha]_\sharp \] where $i(\alpha,\beta)$ is the geometric intersection number.

\end{lem}

The lemma below provides sufficient conditions for a simple closed curve in $S$ to be spiraling and admissible.

If $\gamma$ and $\gamma'$ are oriented admissble curves and $\rho(\gamma) = \rho(\gamma')$, we say that they are oriented \textit{consistently} if their respective developed images have components $\gamma_D$ and $\gamma_D'$ such that with the induced orientation they have the same initial endpoint.

\begin{lem} \label{alt int num}
Assume $\gamma$, $\gamma'$ and $\lambda$ are simple closed curves in a standard projective structure $\Sigma(\rho,\lambda)$ such that $\gamma \cap \lambda = \phi$, $\rho(\gamma) = \rho(\gamma')$, $\gamma$ and $\gamma'$ are oriented consistently and  \[|\hat{i}(\gamma, \gamma')| = i(\gamma,\gamma') = \frac{1}{2} i(\gamma',\lambda). \]  If for every arc $a' \in \gamma' - \gamma$ there is an arc $a \in \gamma - \gamma'$ such that the curve $a' \cup a$ has trivial holonomy, then $\gamma'$ is spiraling and admissible in $\Sigma(\rho,\lambda)$.
\end{lem}

    \noindent \textbf{Proof : } 

Since $\gamma \cap \lambda = \phi$, Lemma \ref{disjoint admissibility} implies  $\gamma$ is admissible.  The standard developing map is a covering map onto its image so $\gamma'$ is admissible since it is simple.  Let $\uc{\gamma}$ and $\uc{\gamma}'$ be components of the universal cover of $\gamma$ and $\gamma'$ which intersect.  

Let $a$  be an arc of $\gamma' - \gamma$.  Since $\lambda$ is the real curve of $\Sigma(\rho,\lambda)$ it is a separating curve so the arc $a$ must cross $\lambda$ an even number of times. 
If $a \cap \lambda$ consisted of more than two points some other arc of $\gamma' - \gamma$ would have to be disjoint from $\lambda$ since $i(\gamma,\gamma') = \frac{1}{2} i(\gamma',\lambda)$.  

Let $\uc{a} \subset \chat$ be the developed image of a component of $\pi^{-1}(a)$ where $\pi : \uc{S} \rightarrow S$ is the universal covering.  Let $D(\uc{\gamma})$ be the developed image of a component of $\pi^{-1}(\gamma)$ such that the arc $\uc{a}$ has both endpoints on $D(\uc{\gamma})$.  Such a $D(\uc{\gamma})$ exists by the assumption that $a' \cup a$ has trivial holonomy.

Assume $D(\uc{\gamma})$ has endpoints $\{0,\infty\}$ and is oriented from $0$ to $\infty$.  Then $D(\uc{\gamma})$ splits the upper half plane into two regions, called region $0$ (which is bounded by $R^+$ and $D(\uc{\gamma}$)) and region $1$   
(which is bounded by $R^{-1}$ and $D(\uc{\gamma})$)).  The lower half plane we call region $2$.  The assumption that $|\hat{i}(\gamma,\gamma')| = i(\gamma,\gamma')$ implies that all indices of $D(\uc{\gamma}) \cap D(\uc{\gamma}')$ are the same.  Assume $\hat{i}(\gamma, \gamma') > 0$.  Then the oriented arc $\uc{a}$ starts in region 1 and ends in region 0 and $\uc{a}$ must pass through region 2 since it cannot intersect $D(\uc{\gamma})$.     This implies that  $\uc{a}$ intersects $\R$ and therefore $a$ must intersect $\lambda$.  By the comment above, since every arc of $\gamma' - \gamma$ intersects $\lambda$ it must do so exactly twice.  The first point of $\uc{a} \cap \R$ is negative and the last point is positive since $\uc{a}$ starts in region 1 and ends in region 0.
Thus $\gamma'$ is right-spiraling.  If $\hat{i}(\gamma, \gamma') < 0$ the argument above implies that $\gamma'$ is left-spiraling.  
\qed

\begin{lem} \label{spiraling twists}
Assume $\gamma \cap \lambda = \phi$ and $\beta$ is a collection of disjoint meridians in $S$
    such that $| \gamma \cap \beta_i | = 1$ and $| \lambda \cap \beta_i| = 2$.  Then 
    $T^{n}_{\beta}(\gamma)$ has a left (or right) admissible spiraling 
    representative in the standard projective structure $\Sigma(\rho, \lambda)$ 
    for $n \in \mathbb{Z^+}$ (or $n \in \mathbb{Z^-}).$  
    \end{lem}
   
     \noindent \textbf{Proof : } 
     
Since $\gamma$ is essential, Lemma \ref{disjoint admissibility} implies $\gamma$  is admissible in 
     $\Sigma(\rho, \lambda).$ 
Let $S$ have the orientation induced by the standard orientation of $\mathbb{C}$.  
Since each $\beta_i$ is a meridian in $S$, $\rho(\beta_i) = 1$ and therefore $\rho(T^n_\beta(\gamma)) = \rho(\gamma)$.  
The fact that $|\gamma \cap \beta_i| = 1$ implies that $T^n_{\beta_i}(\gamma) = [\gamma, n \beta_i]_\flat$.  Since $\hat{i}(\gamma,T^n_{\beta}(\gamma)) = -n$ we have
\[| \hat{i}(   [T^n_\beta(\gamma) , \gamma ) | = i( T^n_\beta(\gamma) , \gamma ) = n |\beta | \]
where $|\beta|$ is the number of components of $\beta$.

\begin{figure}[ht]
	  \psfrag{1}{\small{$p'_{1}$}}
	  \psfrag{2}{\small{$p'_{2}$}}
	  \psfrag{3}{\small{$p'_{3}$}}
	  \psfrag{4}{\small{$p'_{4}$}}
	  \psfrag{5}{\small{$p'_{5}$}}
	  \psfrag{6}{\small{$p'_{6}$}}
	 \centerline{ \epsfig{file=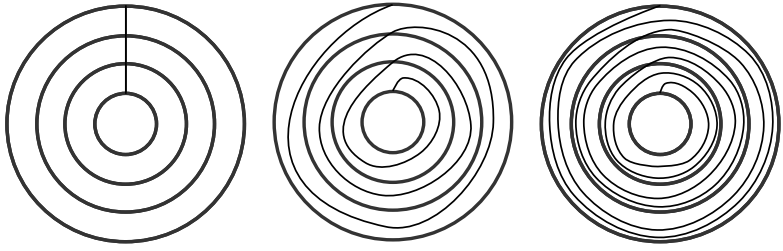,scale=.8} }
	    \caption{ Grafting along left and right alternating curves realized as $\sharp$ and $\flat$ operations.}
	    \label{fig: meridian twists}
	    \end{figure}

Since $| \lambda \cap \beta_i| = 2$ each twist of $\gamma$ about $\beta$ adds $2$ intersection points to $i (T^n_{\beta_i} (\gamma), \lambda)$, so
\[ i ( T^n_{\beta_i} (\gamma), \lambda ) = 2n | \beta | .\]
In \cite{Ito} Ito observed that $[\gamma, T^n_\beta(\gamma)]_\sharp = [\gamma, [\gamma, n \beta]_\flat]_\sharp = n \beta$.  There are exactly $n | \beta |$ arcs of both $T^n_\beta(\gamma) - \gamma$ and $\gamma - T^n_\beta(\gamma)$.  Therefore 
each component of $n \beta$ is composed of exactly $1$ arc of $T^n_\beta(\gamma) - \gamma$ and one arc of  $\gamma - T^n_\beta(\gamma)$.  This implies that for every arc of $T^n_\beta(\gamma) - \gamma$ there is an arc of $\gamma - T^n_\beta$ whose union has trivial holonomy, see Figure \ref{fig: meridian twists}.  Since a Dehn twist is an orientation preserving homeomorphism for any orientation of $\gamma$, then with the induced orientation on $T^n_\beta(\gamma)$, $\gamma$ and $T^n_{\beta}(\gamma)$ are oriented consistently.
   
Now we can apply Lemma \ref{alt int num} to see that $T^n_{\beta}(\gamma)$ is spiraling.  Since  
$\hat{i}(\gamma, T^n_\beta(\gamma)) = -n$, Lemma \ref{alt int num} implies that if $n > 0$, $T^n_{\beta_i}(\gamma)$ is left-spiraling, and if $n > 0$,  $T^n_{\beta_i}(\gamma)$ is right-spiraling.

 \qed

%

Next we use Theorem \ref{flatsharp} with Lemma \ref{twists and sharp} to show that in some cases the real lamination produced by grafting is a Dehn twist of the original lamination along the grafting curve.

\begin{cor} \label{dehn twists}
    Suppose $\gamma$ is an admissible right-spiraling 
    representative for which \\
    $| \lambda \cap \gamma | = 2$. Then    
    \[Gr_{\gamma}(\Sigma(\rho,\lambda)) = \Sigma(\rho,  T_{\gamma}(\lambda))
 .\]  
    If $\gamma$ is an admissible left-spiraling 
    representative for which 
    $| \lambda \cap \gamma | = 2$. Then    
    \[Gr_{\gamma}(\Sigma(\rho,\lambda)) = \Sigma(\rho,  T^{-1}_{\gamma}(\lambda))
 .\]

    \end{cor}
    
    \noindent \textbf{Proof : } 

Theorem \ref{flatsharp} implies the new real curves are given by $\lambda_\gamma \simeq (\lambda, 2 \gamma)_\flat$ in the case that $\gamma$ is right-spiraling and 
$\lambda_\gamma \simeq (\lambda, 2 \gamma)_\sharp$ if $\gamma$ is left-spiraling.  Since $| \lambda \cap \gamma | = 2$, Lemma \ref{twists and sharp} implies $ T_\gamma (\lambda) = (\lambda, 2\gamma)_\flat$ if $\gamma$ is right-spiraling and 
$ T_\gamma (\lambda) = (\lambda, 2\gamma)_\sharp$. if $\gamma$ is left-spiraling.

 \qed

	%
	%
	%
	%
    
\begin{figure}[ht]
\begin{center}
	    \psfrag{2}{$\sharp$}
	    \psfrag{1}{$\flat$}
	 	\includegraphics[scale=.8]{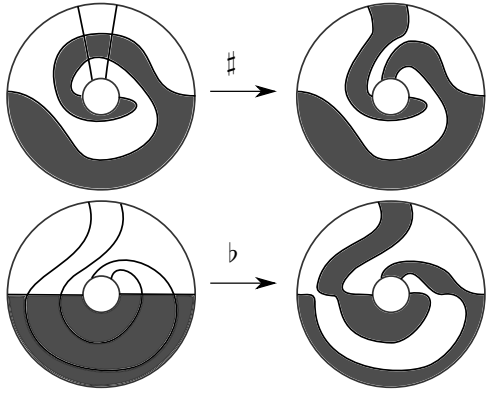}
	    \caption{ The annulus $A$ where $ [\lambda, T^{k}(2\gamma)]_{\sharp} \simeq [T^{k}(\lambda), 2 \gamma]_{\flat}$ for $k = -1$.}
	    \label{figure: twistsandsharps}
\end{center}
	    \end{figure}

Our next goal is to 
show that for 
every two structures
that differ 
by a power of an elementary twist (via $T^{k}_{\beta}$),
there are an infinite number of structures realized as a graft 
of each.  How this is done is suggested in Figure \ref{figure: twistsandsharps}.

The next few computations will make use of the following formula for the geometric intersection number $i(a,b)$ for curves $a$ and $b$ on a torus \cite{FM}.

\[ i(T^k_{(p,q)}(r,s),(r,s)) = |k| \cdot |ps - qr|^2 .\]

\begin{lem} \label{sharp = flat} Let $\beta, \gamma, 
    \lambda$ be simple closed multi-curves in $S$. 
    Assume $| \beta \cap \lambda | = 2,  |\beta \cap \gamma| = 1$ and $|\gamma \cap \lambda| = 0$.  
	Set $T^{k} = T^k_{\beta_T}$ for $k \in \mathbb{Z}$.  

Then
	    \[ [\lambda, T^{k}(2\gamma)]_{\sharp} \simeq
	    [T^{k}(\lambda), 2 \gamma]_{\flat} \] and
	    \[ T^{k}[\lambda, 2 \gamma]_{\sharp} \simeq [\lambda,
	    T^{2k}(2\gamma)]_{\sharp} \simeq [T^{2k}(\lambda), 2
	    \gamma]_{\flat} \simeq [T^{k}(\lambda), T^{k}(2
	    \gamma)]_{\sharp} .\]
	    \end{lem}
	     
\noindent \textbf{Proof : } 
%

  Since $|\gamma \cap \lambda| = 0$ there is an element $b \in [\beta]$ and a regular neighborhood $A$ of $b$ containing all
  intersection points $\lambda \cap T^k (\gamma)$.   After replacing $\beta$ with $b$ we may assume that
  $(\lambda \cap T^{2k}(\gamma))$ , $(\lambda \cap T^{k}(\gamma))$ and
  $(T^{k}(\lambda) \cap \gamma)$ are all contained in $A$.  Then for $k
  \in \mathbb{Z}$ we have \[ [\lambda, T^{k}(2\gamma)]_{\sharp}
  \arrowvert_{(S \backslash A)} \simeq [T^{k}(\lambda), 2 \gamma]_{\flat}
  \arrowvert_{(S \backslash A)} \simeq ( \lambda \cup 2 \gamma ) \arrowvert_{(S
  \backslash A)} \] and 
  \[ [T^{2k}(\lambda), 2 \gamma]_{\flat} \arrowvert_{(S \backslash A)} \simeq T^{k}[\lambda, 2 \gamma]_{\sharp} \arrowvert_{(S \backslash A)}
  .\]
		
 For a curve $\alpha \in S$ let $\alpha_T = \psi(\alpha \cap A)$ where $\psi : A \rightarrow T$ is a map of $A$ similar to the the map in Lemma \ref{sharp passes} which 
identifies $\partial A$ so that $\psi(\lambda \cap A) = (2,0)$, $\psi(\beta \cap A) = (0,1)$ and $\psi(\gamma \cap A) = (1,0)$. 
Then $T^k_{\beta_T}(\gamma_T) = (1,k)$ and $T^k_{\beta_T}(\lambda_T) = (2,2k)$.  
 
 Now \[ [\lambda_T, T^k_{\beta_T}(2 \gamma_T)]_\sharp  = [(2,0), (2,2k)]_\sharp = (4,2k) .\] 
 And since \[ [T^k_{\beta_T}(2 \gamma_T), \lambda_T]_\sharp = [(2,2k), (2,0)]_\sharp = (4,2k)\] 
Lemma \ref{flat sharp switch} implies that \[ [\lambda_T, T^k_{\beta_T}(2 \gamma_T)]_\sharp = [T^k_{\beta_T}(2 \gamma_T), \lambda_T]_\flat .\]

	Also
	\[ T^{k}[\lambda_{T}, 2 \gamma_{T}]_{\sharp} =
	T^{k}[(2,0),(2,0)]_{\sharp} =
	T^{k}(4,0) = (4,4k) \]
	and
	\[ [\lambda_{T}, T^{2k}(2 \gamma_{T})]_{\sharp} = 	
	    [(2,0),(2,4k)]_{\sharp} = (4, 4k) \] 
	    and
	\[ [T^{2k}(\lambda_{T}), 2 \gamma_{T}]_{\flat} =
	    [(2,0),(2,4k)]_{\sharp} = (4, 4k) \]
              and
         	\[ [T^k (\lambda_{T}), T^{k}(2 \gamma_{T})]_{\sharp} = 	
	   [(2,2k),(2,2k)]_{\sharp} = (4, 4k) \]
 	 and
            so by Lemma \ref{sharp passes} \[
	    [\lambda, T^{2k}(2\gamma)]_{\sharp} 
	    	    	\simeq [T^{2k}(\lambda), 2 \gamma ]_{\flat}
			 \simeq T^{k}[ \lambda,   2\gamma]_{\sharp} 
			 \simeq [ T^k (\lambda), T^k (2 \gamma)]_{\sharp} \] 
		      \qed

   

  The final theorem provides a sufficient condition on an admissible curve $\delta$
    for the real  lamination obtained by grafting $\Sigma$ along $\delta$ to be the
    multicurve $T^{k}_{i}(\lambda) \cup 2\delta$.         
 Furthermore, it asserts that for any two structures 
$\Sigma(\rho, \lambda)$ and $\Sigma(\rho, T^{k}_{i}(\lambda))$ 
differing by an iterated elementary 
twist there are an infinite number of ways (parameterized by $m \in 
\mathbb{Z}$)
to connect them by grafting.  That is, for each $k$ there are a pair of admissible spiraling curves such that grafting $\Sigma(\rho, \lambda)$ along one and $\Sigma(\rho, T^{k}_{i}(\lambda))$ along the other produce the same projective structure.  Therefore, there are an infinite 
number of structures realized as a graft of both
$\Sigma(\rho, \lambda)$ and $\Sigma(\rho, T^{k}_{i}(\lambda))$.

%

\begin{thm} \label{iterated connectivity}
    Let $k,l,m$ be integers such that $k+l = m$.  Let $\beta, \gamma, 
    \lambda$ be simple closed multi-curves in $S$. 
    Assume $\gamma$ is an admissible representative of $[\gamma]$ in $\Sigma(\rho, 
    \lambda)$ and $| \beta \cap \gamma | = 1$, $| \beta \cap \lambda 
    | = 2$, and $\lambda \cap \gamma = \phi$.  Then    
    \[ Gr_{T^{2m}_{\beta}(\gamma)}(\Sigma(\rho, \lambda)) = 
    Gr_{T^{k}_{\beta}(\gamma)}(\Sigma(\rho, T^{l}_{\beta}(\lambda)) . 
    \]
    \end{thm}
    
    \noindent \textbf{Proof : } 
	Assume $m > 0$.  Lemma \ref{spiraling twists}
	implies $T^m_\beta(\gamma)$ is a right-spiraling  admissible representative in  
	$\Sigma(\rho, \lambda)$.  Then $T^k(\gamma)$ is right-spiraling in
	$\Sigma(\rho, T^{-l}_{\beta}(\lambda))$ and this implies that  
	$T^k(\gamma)$ has a left-spiraling admissible representative in
$\Sigma(\rho, T^{l}_{\beta}(\lambda))$.

	Let $Gr_{T^{2m}_{\beta}(\gamma)}(\Sigma(\rho, \lambda)) = 
	\Sigma(\rho, \lambda_{1})$ and $
	Gr_{T^{k}_{\beta}(\gamma)}
	(\Sigma(\rho, T^{l}_{\beta}(\lambda)) = \Sigma(\rho, 
	\lambda_{2})$.  Theorem \ref{flatsharp} implies 
	\[ \lambda_{1} = [\lambda, 2 T^{2m}_{\beta}(\gamma)]_{\flat} 
	\mbox{ and } \lambda_{2} = [T^{l}_{\beta}(\lambda), 2 
	T^{k}_{\beta}(\gamma)]_{\sharp} .\]
	Then Lemma \ref{sharp = flat} implies \[ \lambda_{2} = 
	[\lambda, 2 T^{k+l}_{\beta}(\gamma)]_{\flat} .\] 

	 A symmetric proof works for $m < 0$.  
\qed


\begin{cor} \label{cor: infinite connectivity}
Any two Schottky projective structures differing by an iterated elementary twist can be connected in an infinite number of ways by grafting.
\end{cor}

\noindent \textbf{Proof : }
Assume $\Sigma(\rho,\lambda_1)$ and $\Sigma(\rho,\lambda_2)$ are projective structures on $S$ and $\lambda_1 = T^k_\beta(\lambda_2)$.  Then for every $m \in \mathbb{Z}$ there is an $l \in \mathbb{Z}$ such that $m = k+l$.  
By Theorem \ref{iterated connectivity} each pair $(m,l)$ gives a different way to connect the two structures by grafting. \qed

\begin{cor} \label{cor: elementary connectivity}
There are an infinite number of standard projective structures that can be connected by grafting.
\end{cor}

\noindent \textbf{Proof : }
Given a standard projective structre  $\Sigma(\rho,\lambda)$,  let $\beta$ be a meridian of $S$.  Since $\beta$ develops to a closed curve in $\chat$ twists about $\beta$ correspond to twisting in $\chat$ as in Figure \ref{fig: meridian twists}.  So for every $l \in \mathbb{Z}$ the structure  $\Sigma(\rho,T^l_\beta(\lambda))$ is a standard projective structure.  
By Theorem \ref{iterated connectivity} these two structures can connected by grafting. \qed

The last corollary relates grafting and twists about two different meridians.
Let $T^i_j = T^i_{\beta_j}$ where $T^{i}$ is $i$ Dehn twists 
about $\beta_{j}$, a meridian in $S$.

\begin{cor} \label{iterated connectivity cor}
    Let $k,l,m$ be integers such that $k+l = m$ and $l \cdot k \geq 
    0$ and $l < k$.  Let $\beta_1, \beta_2 \gamma, 
    \lambda$ be simple closed multi-curves in $S$ and assume $\gamma$ is admissible in $\Sigma(\rho, \lambda)$.
  Assume $| \beta_i \cap \gamma | = 1$, $| \beta_i \cap \lambda 
    | = 2$, and $\lambda \cap \gamma = \phi$.
    Then \[ Gr_{T^{k}_{1} T^{l}_{2} (\gamma) } (\Sigma(\rho, 
    \lambda)) =  Gr_{T^{l}_{2}(\gamma)} (\Sigma(\rho, 
    T^{k}_{1}(\lambda))).
    \]
    \end{cor}
    
    \noindent \textbf{Proof : } 
     Assume first that $l$ and $k$ are both positive.  By Lemma \ref{spiraling twists} 
     $T^{k}_{1} T^{l}_{2} (\gamma)$ has an admissible right-spiraling 
    representative in $\Sigma(\rho, \lambda)$ 
Since $l<k$ it also follows from Lemma \ref{spiraling twists} that 
     $T^{l}_{2}(\gamma)$ is left-spiraling relative to 
     $T^{k}_{1}(\lambda)$. Then 
     \[ Gr_{T^{k}_{1} T^{l}_{2} (\gamma) } (\Sigma(\rho, 
    \lambda)) = Gr_{T^{l}_{2}(\gamma)} (\Sigma(\rho, 
    T^{k}_{1}(\lambda)))  
     \] by Theorem \ref{iterated connectivity}. 
     
    If $l \cdot k = 0$ and $k = 0$ then $T^{k}_{1} (\delta) =
    \delta$ for all closed curves $\delta$ in $S$.  Then $T^{k}_{1} 
    T^{l}_{2}(\gamma) = T^{l}_{2}(\gamma)$ is left-spiraling by
    Lemma \ref{spiraling twists}. Then 
    \[  Gr_{T^{k}_{1} T^{l}_{2} (\gamma) } (\Sigma(\rho, 
    \lambda)) = Gr_{T^{l}_{2}(\gamma)} (\Sigma(\rho, \lambda))  
    \] 
by 
    Theorem \ref{iterated connectivity}.  A similar proof works for 
    the case that $l =0$. \qed



\begin{thebibliography}{1}
	    
		\bibitem[1]{Baba} S. Baba, \textit{Complex Projective Structures with Schottky Holonomy}, Geometric And Functional Analysis
Volume 22, Number 2 (2012), 267-310, DOI: 10.1007/s00039-012-0155-x.

	
	     \bibitem[2]{Bon} F. Bonahon, \textit{Geodesic laminations on surfaces}, 
	 	Contemp.  Math. \textbf{269} (2001) 1-37.
	     
	     \bibitem[3]{CB} A. J. Casson \& S. A. Bleiler, \textit{Automorphisms  of 
	     surfaces after Nielsen and Thurston}, London Mathematical Socieity 
	     Student Texts No. 9, Cambridge University Press, Cambridge, 1988.
	 
	 
	     \bibitem[4]{Ear} Clifford Earle, \textit{On variation of projective structures}, Riemann 
	     surfaces and related topics, Ann. Math. Studies \textbf{97} (1981) 87-99.
	     
	     \bibitem[5]{FM}B.Farb \& D.Margalit, \underline{A primer on mapping class groups}, Princeton Mathematical Series, Princeton University Press, 2011.
	     
	     \bibitem[6]{GKM} D. Gallo, M. Kapovich \& and A. Marden,  \textit{The
	     monodromy groups of Schwarzian equations on closed Reimann surface},
	     Ann. of Math \textbf{151} (2000) 625 - 704 
	 
	     \bibitem[7]{Gol} W. Goldman, \textit{Projective structures with Fuchsian
	     holonomy}, J. Diff. Geom \textbf{25} (1987) 297 - 326.
	 
	 	\bibitem[8]{Pap} U. Hamenst\"{a}dt, \textit{Geometry of the complex of curves and of Teichm\"{u}ller space}, Handbook of Teichm\"{u}ller theory Vol. 2,  European Mathematical Society, Zurich 2007. 	
	 %
	     \bibitem[9]{Hat} A. Hatcher, \textit{Algebraic Topology}, Cambridge University Press, Cambridge, 2002.
	 
	     \bibitem[10]{Hej} D. Hejhal, \textit{Monodromy Groups and Linearly 
	     Polymorphic Functions}, Acta Math \textbf{115} 
	     (1975) 1-55.
	     
	     \bibitem[11]{Hub} J.H. Hubbard, \textit{The monodromy of projective 
	     structures}, Riemann surfaces and related topics, Ann. 
	     Math. Studies \textbf{97} (1981) 257-275.
	     
	     \bibitem[12]{Ito} K. Ito, \textit{Exotic projective structures and
	     quasifuchsian spaces II}, Duke Math J. \textbf{140} No. 1 (2007) 85-109.
	 
	 	
	     \bibitem[13]{KT} Y. Kamishima \& S.P. Tan, \textit{Deformation spaces on 
	     geometric structures}, Advanced Studies in Pure Mathematics
	     \textbf{20} (1992) 263-299.
	 

	     \bibitem[14]{Kap89} M. Kapovich, \textit{On conformal structures with Fuchsian holonomy},
		 Soviet Math. Dokl. Vol. \textbf{38} (1989) N 1 14-17. 
		

	
	 	\bibitem[15]{Kap} M. Kapovich, \textit{Hyperbolic manifolds and discrete groups},  Progress in Mathematics \textbf{183} (2001) 1-468. 
	     
	         \bibitem[16]{MT} K. Matsuzaki \& M. Taniguchi, \textit{Hyperbolic manifolds 
	     and Kleinian groups}, Oxford University Press 1998.
	    
		\bibitem[17]{Tan88} S. Tan, \textit{Representations of surface groups into $PSL_2(\R)$ and geometric structures}. PhD thesis, University of California, Los Angeles, 1988.

	     \bibitem[18]{Thu} W. Thurston,  \textit{Geometry and topology of 3-manifolds}, Vol. 1, Princeton University Press, Princeton, 1977.
	
	 \end{thebibliography}



\end{document}